\documentclass[11pt,notitlepage,twoside]{article}
\usepackage{color}

\usepackage{latexsym}
\usepackage{amsmath}
\usepackage{amsfonts}
\usepackage{amssymb}
\usepackage{enumerate}
\usepackage[latin1]{inputenc} 
\renewcommand{\a }{\alpha }

\renewcommand{\d}{\delta }

\newcommand{\D }{\Delta }

\newcommand{\e }{\varepsilon }
\newcommand{\eps }{\varepsilon }
\newcommand{\g }{\gamma}

\renewcommand{\l }{\lambda }

\newcommand{\n }{\nabla }

\renewcommand{\o }{\omega }
\renewcommand{\O }{\Omega }

\newcommand{\wtilde }{\widetilde}

\newcommand{\be}{\begin{equation}}
\newcommand{\ee}{\end{equation}}
\newenvironment{pf}{\noindent{\bf Proof.}\enspace}{
\hfill$\Box$\medskip}

\newcommand{\R}{\mathbb{R}}
 
\newtheorem{thm}{Theorem}[section]
\newtheorem{pro}[thm]{Proposition}
\newtheorem{lem}[thm]{Lemma}
\newtheorem{rem}[thm]{Remark}

\numberwithin{equation}{section} 

\title { \Large \textbf{ Non existence of  solutions for a slightly super-critical \\
elliptic problem with non-power nonlinearity }}
\author{{\bf\large Mohamed Ben Ayed}
{\bf\large}\vspace{1mm}\\
 {\it\small Department of Mathematics,  College of Science, Qassim University, Saudi Arabia,}\\
{\it\small e-mail: M.Benayed@qu.edu.sa}\\
 {\it\small Universit\'e de Sfax, Facult\'e des Sciences de Sfax, Route Sokra, Sfax, Tunisie, } \\
{\bf\large Habib Fourti} 
{\bf\large}\vspace{1mm}\\
{\it\small Department of Mathematics and Statistics, College of
Science,}\\
{\it\small King Faisal University, Al-Ahsa 31982, Saudi Arabia,}\\
{\it\small e-mail: hfourti@kfu.edu.sa
} }

\begin{document}

\date{}

\maketitle

{\footnotesize \noindent {\bf Abstract.} In this paper, we are
concerned with the following elliptic equation
$$ ( SC_\e ) \qquad
 \begin{cases}  -\Delta u =  |u|^{4/(n-2)}u [\ln (e+|u|)]^\e & \hbox{ in } \Omega,\\
u = 0 & \hbox{ on }\partial \Omega,
\end{cases}
$$ where $\Omega $
 is a smooth bounded open domain in $\mathbb{R}^n, \ n\geq 3$ and
 $\e>0$. In Comm. Contemp. Math. (2003), Ben Ayed et al. showed that the slightly
 supercritical usual elliptic problem has no single peaked solution.
 Here we extend their result for problem $( SC_\e )$ when $\e$ is small
 enough, and that by assuming a new assumption.

 }


{\footnotesize\section{Introduction and results}


 Let $\Omega $ be a smooth bounded open domain in $\mathbb{R}^n, \ n\geq 3$ and
 $\e>0$. We consider the following nonlinear elliptic problem
$$ ( SC_\e ) \qquad
\begin{cases} \displaystyle - \Delta u =  |u|^{4/(n-2)}u [\ln (e+|u|)]^\e &\hbox{ in } \Omega,\\
u = 0 & \hbox{ on }\partial \Omega . \end{cases}
$$

 This problem is considered slightly supercritical.
Indeed, the critical scenario occurs when $\e = 0$, leading to the
nonlinearity $|u|^{ 4 / (n-2)}u$. In this context, the associated
Euler functional fails to meet the Palais-Smale condition, resulting
the lack of compactness in the associated variational problem. The
problem $( SC_0 )$ that arises, is particularly significant in the
field of geometry, as it relates to the Yamabe problem, which is a
variant of this issue on manifolds. Below, we summarize some
established findings regarding the critical case i.e.  $\e = 0$. In
the case of a star-shaped domain $\O$, Pohozaev demonstrated in
\cite{P} that problem $(SC_0)$ does not yield any positive
solutions. For an annular domain $\O$, Kazdan and Warner showed in
\cite{KW} that a positive radial solution exists. Utilizing the
theory of critical points at infinity, Bahri and Coron \cite{BC}
established that this problem has a
positive solution, provided that $\O$ possesses nontrivial topology.\\
As far as we know problem $(SC\e)$ has not been studied in
literature yet. In this paper our main goal is to establish the
nonexistence of single peaked solution for our problem. This is in
contrast with the slightly subcritical scenario when studying the
problem
$$ \displaystyle ( P_\e ) \qquad
\begin{cases}\displaystyle  - \Delta u =  \frac{|u|^{4/(n-2)  }u}{ [\ln (e+|u|)]^\e }  &\hbox{ in } \Omega,\\
u = 0 & \hbox{ on }\partial \Omega . \end{cases}
$$

 To present both established and novel findings, it is
beneficial to revisit some familiar definitions.\\
 The space $H^1_0(\O)$ is equipped with the norm $\|.\|$ and its
corresponding inner product $\left<.,.\right>$ defined by
$$ \|u\|^2=\int_{\O}|\n
u|^2;\quad \left<u,v\right>=\int_{\O}\n u \n v,\quad u,v \in
H^1_0(\O).$$ For $a\in \O$ and $\l>0$, let
 \begin{equation}\label{d}
\d_{ a,\l }(y)=\frac{c_0\l^{(n-2)/2}}{\left(1+\l^2|y-a|^2\right)^{(n-2)/2}},
\hbox{ where } c_0:=(n(n-2))^{(n-2)/4}.
 \end{equation}
The constant $c_0$ is chosen such that $\d_{(a,\l)}$ is the family
of solutions of the following problem
\begin{equation}\label{e}
 -\Delta u=u^{(n+2)/(n-2)},\;\;u>0\; \;\mbox{in}\; \;\R^n.
\end{equation}
Notice that the family $ \d_{(a,\l)}$ achieves the best Sobolev
constant \begin{equation}\label{S_n}S_n:=\inf\{\|\n
u\|^2_{L^2(\R^n)}\|u\|^{-2}_{L^{2n/(n-2)}(\R^n)}:u\not\equiv 0, \n u
\in (L^{2}(\R^n))^n \,\mbox{and}\;u \in L^{2n/(n-2)}(\R^n)\}.
\end{equation}
We denote by $P\d_{(a,\l)}$ the projection of $\d_{(a,\l)}$ onto $
H^1_0(\O)$, defined by
\begin{eqnarray}\label{q1}
-\Delta P\d_{ a,\l }=-\D \d_{ a,\l } \mbox{ in } \O, \quad
P\d_{ a,\l }=0 \mbox{ on}\;
\partial \O.
\end{eqnarray}
We will denote by $G$ the Green's function and by $H$ its regular
part, that is
$$G(x,y)=|x-y|^{2-n}-H(x,y) \quad\mbox{ for }\,\,
(x,y)\in \O^2,$$ and for $x\in\O$, $H$ satisfies
\begin{eqnarray*}   \begin{cases}
\Delta  H(x,.)=0   & \mbox{in} \, \, \O,
\\H(x,y)=|x-y|^{2-n},&\mbox{for}\, \,y\in
\partial \O.
\end{cases}
\end{eqnarray*}

We define the Robin function as $$ R(x) = H(x, x),\ x\in \O . $$

 The initial study demonstrating the presence of blowing-up solutions to
problem $(P_\e)$ is referenced as \cite{Pistoianew}. In this work,
the authors established that any non-degenerate critical point $x_0$
of the Robin function $R$ gives rise to a set of single-peak solutions that concentrate
around $x_0$ as $\e$ approaches $0$.  This family of solutions can
be expressed in the following manner:
$$u_\e = \a_\e P\d_{a_\e , \l_\e } + v _ \e,$$
where $ \a_\e \to 1$, $a_\e\in \O
 $, $a_\e \rightarrow x_0$,  $\l_\e \rightarrow \infty$, $v\in E_{(a,\l)} $ and $v_\e\rightarrow 0 $ in $H^1_0(\O),$ as $\e\rightarrow
 0$,
where
\begin{align}\label{vorthogonal}&E_{(a,\l)}:=
\Big\{v\in H^1_0(\O):\bigg<v,P\d_{a,\l}\bigg>=\left<v,\frac{\partial
P\d_{a, \l}}{\partial \l }\right>= \left<v,\frac{\partial P\d_{a,
\l}}{\partial a_j}\right>=0\; \, \, \forall \; 1\leq j\leq n \Big\},
\end{align}
and where  $ a _j $ is the $j^{th}$ component of $ a $.\\
 Recently, we proved in \cite{BFG} the existence of both
positive solutions and those with changing signs that exhibit
blow-up and/or blow-down behavior at various locations in $\O$. The
solutions we found have the following expansion
\begin{equation*}\label{u_eps}
u_\e= \sum_{i=1}^m\alpha_i\gamma_iP\delta_{ a_i,\lambda_i }+v,
\end{equation*}
where $m$ is a positive integer,
$(\gamma_1,\ldots,\gamma_m)\in\{-1,1\}^m$, $(\a_1,...,\a_m)\in
(0,+\infty)^m$, $(\l_1,...,\l_m)\in (0,+\infty)^m$ and
$(a_1,...,a_m)\in \O^{m}$.  Problem $(P_\e)$ also admits
bubble tower solutions with alternating signs, as shown in
\cite{FG}. These towers require the spikes to concentrate at the
same point, namely at a stable critical point of the Robin
function.\\
Regarding the investigation of the profile of subcritical solutions,
the authors in \cite{Pardo} studied the asymptotic behavior of
radially symmetric solutions of $(P_\e)$ when the domain is a ball.
This analysis was very recently extended in \cite{F} to the case of
a general domain. In fact, the author analyzed the asymptotic
behavior of the least energy solutions of problem $(P_ \e)$ as
$\e\rightarrow 0$. He demonstrated that this family of solutions
blows up at a critical point of
the Robin function.

In this paper, we prove the nonexistence of solution concentrating
at one point for problem $(SC_\e)$. Precisely, our first result can
be stated as follows.
 \begin{thm} \label{1.4} There exists $ \e_0 > 0 $ such that, for each $ \e \in ( 0 , \e_0) $, Problem $ (SC_\e) $ has no solution
 $ u_\e $ which converges weakly to $0$ (not strongly) and   satisfies:
\begin{equation}\lim _{\e \to 0 } \int_\O | u_\e |^{2n / (n-2) } \ln ^{\e n /2}  ( e + | u_\e | ) = S_n^{n/2} .\tag{\textbf{A}}
\end{equation}
\end{thm}
 Theorem \ref{1.4} extend the nonexistence of single peaked solution concerning the usual slightly supercritical elliptic problem with exponent nonlinearity to
 nonpower one. Infact, the authors in \cite{BEGR} were interested in the following problem
\begin{eqnarray}\label{suppb}\begin{cases}
 -\Delta u= | u |^{p-1+\e} u     & \mbox{in} \, \, \O,  \\
 \quad \ \ u= 0     &  \mbox{on}\, \, \partial \O,
\end{cases}
 \qquad \mbox{where } \quad p+1:= 2n / (n-2) \quad \mbox{  and } \quad \e > 0 .
\end{eqnarray}

In the following, we state their result.
\begin{thm}\label{thBEGR}\cite{BEGR} Let $\O$ be any smooth bounded domain in $\mathbb{R}^n, \ n\geq
3$. $\eqref{suppb}$ has no positive solution $u_\e$ such that
 \begin{equation*} u_\e =  P\d_{a_\e , \l_\e } + v _ \e
 \end{equation*} with $v_\e\rightarrow 0 $ in $H^1_0(\O),$ $a_\e\in \O
 $ and $\l_\e d(a_\e,\partial \O)\rightarrow \infty$ as $\e\rightarrow
 0$.
\end{thm}
Compared with Theorem \ref{thBEGR}, our result presents a condition
of new type. Moreover, our theorem is more practicable in the sense
that, checking Assumption $(A)$ is easier than the one of the
mentioned result.\\
We also mention that by following the same argument developed in our
paper, one can prove the following result.
\begin{thm}  There exists $ \e_0 > 0 $ such that, for each $ \e \in ( 0 , \e_0) $, Problem \eqref{suppb} has no solution
 $ u_\e $ which converges weakly to $0$ and   satisfies:
\begin{equation}\lim _{\e \to 0 } \int_\O | u_\e |^{\frac{2n }{n-2 }+ \e \frac{n }{2} }  = S_n^{n/2} .
\end{equation}
\end{thm}
There is a wide literature about the study of the blow up solutions
for the almost critical problem \eqref{suppb} (see for instance
\cite{DFM1,DFM2, KR, Passaseo1,Passaseo2, PR}).
 In references \cite{DFM1} and \cite{DFM2}, Del Pino, Felmer and Musso (also see
Khenissy-Rey \cite{KR} for dimension $3$) developed solutions to
equation \eqref{suppb} blowing up at two points as the parameter
approaches zero. Later, Pistoia and Rey investigated the same
problem \eqref{suppb} in pierced domains. They improved the previous
results by establishing the existence k-peaked solutions (for any
$k\geq 2$) in the case of particular symmetric domain with small
holes and solutions
blowing up at two or three points in the nonsymmetric case provided that $\e$ is small enough.\\
We mention that there is an interesting analogy between the results
obtained for the almost critical problem $(SC_\e)$ and those known
for the usual almost critical elliptic problem \eqref{suppb} in both
subcritical and supercritical regimes. Thus one may expect the
existence of solution of $(SC_\e)$ concentrating at two distinct
points.\\
Recently, the authors in \cite{Benayed,BO} were interested to the
counterpart of \eqref{suppb} when the nonlinearity in \eqref{suppb}
is replaced by $K u^{p+\e}$. Thanks to the  finite-dimensional
reduction of a supercritical exponent equation developed in
\cite{Benayed}, the authors in \cite{BO} were able to construct
solutions by following the ideas of Bahri-Li-Rey \cite{BLR}. We
believe that the same phenomenon holds for the nonpower nonlinearity
$Ku^p\ln(e+u)^\e$. Verification of this conjecture remains as the
future work \cite{BF}.\\
To prove Theorem \ref{1.4}, we argue by contradiction, assuming that a solution $u_\e$ exists. Under the assumption $(A)$ and the fact $ u_\e  $ converges weakly to $0$, we derive that $  \max | u_\e | \to \infty$.  Notice that $ - u_\e $ is also a solution of $(SC_\e)$ and therefore, without loss of the generality, we can assume that $ \max u_\e \to \infty$.  Furthermore, we give the expansion of the family $u_\e$.  More precisely, we get in Proposition \ref{pro2.2} that
$$ u_\e=\a_\e P\d_{a_\e , \l_\e } + v _ \e.$$
Thanks again to assumption $(A)$, we were able to estimate the remainder part $v_\e$. Testing the equation of $(SC_\e)$ by $\lambda_{\varepsilon}\left(\partial P \delta_{a_\e, \l_\e}\right) /\left(\partial \lambda_{\varepsilon}\right)$, we get a balancing condition for the  concentration parameters of $u_\e$ which leads to a contradiction. We point out that the existence of the non power nonlinearity with the critical exponent will lead to some difficulties. To overcome this issue, some quite involving estimates are needed.

Our paper is organized as follows. The next section will be devoted
to some preliminary results. In Section $3$, we study the $v$-part
of a possible solution. Lastly, Theorem \ref{1.4} is proved in
Section $4$.
\section{Preliminary results}
Let $ u_\e$ be a solution of $( SC_\e )$ satisfying Assumption
$\textbf{(A)}$.  Notice that, multiplying the first equation of  $
(SC_\e) $ by $ u_\e $ and integrating over $ \O $, it holds

 \begin{equation}\label{eq01}
   \int_\O | \n u_\e |^2 = \int_\O | u_\e |^{2n/(n-2) } \ln ^\e ( e + | u_\e | ) .
  \end{equation}

We start with the following result.
\begin{lem}\label{lemma 3.1} Let $u_{\varepsilon}$ satisfy the assumption of Theorem \ref{1.4}. Then we have
$$  \int_\O | u_\e |^{2n/(n-2) }  \leq S_n ^{n/2} + o(1)  ,  \leqno{(1)}   $$
$$  \|  u_\e \|  = ( \sqrt{S_n} + o(1) )  \| u_\e \|_{ L^{2n/(n-2) } }  ,  \leqno{(2)}  $$
 as $\varepsilon \rightarrow 0$, where $ S_n$ denotes the Sobolev constant defined by \eqref{S_n}.
 \end{lem}
 \begin{pf}
Using $(A)$ and since $ \ln ( e + | u_\e | ) \geq 1$, it follows
that
 $$ \int_\O | u_\e |^{2n/(n-2) }  \leq \int_\O | u_\e |^{2n/(n-2) }  \ln ^{\e n/2 }  ( e + | u_\e | )  \leq S_n^{n/2} + o(1) $$
which implies the first assertion. \\For the second assertion,
\eqref{eq01} and Holder's inequality together with assumption
$\textbf{(A)}$ assert that
\begin{align*}  \int _\O | \n u_\e |^2 =   \int_\O | u_\e |^{2n/(n-2) }  \ln ^{\e }  ( e + | u_\e | ) & = \int_\O | u_\e |^2 | u_\e |^{4/(n-2) }  \ln ^{\e }  ( e + | u_\e | ) \\
 & \leq \Big( \int_\O | u_\e |^{2n/(n-2) } \Big)^{(n-2)/n} \Big( \int_\O | u_\e |^{2n/(n-2) } \ln ^{\e n /2}  ( e + | u_\e | ) \Big)^{ 2 / n } \\
  & \leq ( S_n ^{n/2} + o(1) )^{2/n} \Big( \int_\O | u_\e |^{2n/(n-2) } \Big)^{(n-2)/n}.
\end{align*}
Thus, by the Sobolev embedding theorem we get
\begin{equation}\label{Quotient}
S_n \leq \frac{ \int _\O | \n u_\e |^2 }{\Big(
\int_\O | u_\e |^{2n/(n-2) } \Big) ^{ (n-2) / n } } \leq S_n + o(1)
.
\end{equation}
 Hence, the second claim follows.
\end{pf}

\begin{pro} \label{minim}
Let $$ I_0 (u) := \frac{1}{2} \int_\O | \n u |^2 - \frac{ n-2 }{2n}
\int_\O | u | ^{2n / (n-2) } \qquad u \in H_0^1 ( \O) .$$ and let
$(u_k)$ be a  sequence satisfying
\begin{equation}\label{integralbound}\lim _{k\to \infty } \| u_k \|^2 = S_n ^{n/2} \qquad \mbox{ and } \qquad \lim _{k\to \infty } \int_\O | u_k | ^{2n / (n-2) }  = S_n ^{n/2} .
\end{equation}
Then, it follows that
$$  \lim _{k\to \infty } I_0 (u_k) = \frac{1}{n} S_n^{n/2} \qquad \mbox{ and } \qquad \lim _{k\to \infty } \| \n I_0(u_k) \| = 0 . $$
\end{pro}
\begin{pf}
First, it is easy to see that if $(u_k)$ satisfies
\eqref{integralbound} then $$\lim _{k\to \infty }  I_0 (u_k) =
\frac{1}{n} S_n^{n/2}. $$ Now, observe that
$$ \n I_0 (u) = u - ( - \D^{-1} ) ( | u | ^{ 4 / (n-2) } u ) $$
where $ ( - \D^{-1} ) ( f ) $ denotes the unique  solution of the
following problem
$$ - \D w = f \quad \mbox{ in } \O \qquad \mbox{ and } \qquad w = 0 \quad \mbox{ on } \partial \O . $$
Thus, let
$$ \o _k := ( - \D^{-1} ) ( | u_k | ^{ 4 / (n-2) } u_k ), $$ it follows that
\be \label{eq2} \| \n I_0(u_k) \| ^2 = \| u_k - \o_k \|^2 = \| u_k
\|^2 + \|  \o_k \|^2 - 2  \langle u_k , \o_k \rangle . \ee Note
that, by \eqref{integralbound}, we get
\begin{equation}\label{eq3} \langle u_k , \o_k \rangle = \int _\O u_k ( - \D \o_k) = \int_\O | u_k |^{2n / (n-2) } \to  S_n ^{n/2} \quad ( \mbox{as } k \to \infty ).
\end{equation}
Using Holder inequality and Sobolev embedding theorem with
\eqref{integralbound}, we have
\begin{align*}
 \|  \o_k \|^2 = \int _\O \o_k ( - \D \o_k) =  \int _\O \o_k  | u_k | ^{ 4 / (n-2) } u_k & \leq \Big( \int _\O | \o_k | ^{2n / (n-2)} \Big)^{(n-2)/(2n)} \Big( \int_\O  | u_k | ^{ 2n / (n-2) } \Big) ^{(n+2)/(2n)} \\
  & \leq \frac{1}{ \sqrt{ S_n} } \| \o_k \| ( S_n ^{n/2} + o(1) )^{(n+2)/(2n)}.
\end{align*}
Thus, it follows that
\begin{equation}\label{eq4}
 \| \o_k \| \leq \frac{1}{ \sqrt{ S_n} } ( S_n ^{n/2} + o(1) )^{(n+2)/(2n)} = S_n ^{n/4} + o(1) .
 \end{equation}
After passing to the limit as $k\rightarrow \infty$, substituting
\eqref{integralbound}, \eqref{eq3} and \eqref{eq4} into \eqref{eq2},
we derive that
$$ \| \n I_0(u_k) \| ^2  = o(1)$$ which completes the proof of the
proposition.
\end{pf}

\begin{pro} \label{pro2.2}Let $u_{\varepsilon}$ be a solution of $\left(SC_{\e}\right)$ which converges weakly to $0$ (not strongly) and satisfies $\textbf{(A)}$. Then   $ u _\e$ has to be written as
 \begin{equation}\label{3.0} u_\e =  \gamma  \a_\e P\d_{a_\e , \l_\e } + v _ \e
 \end{equation} with  $\gamma \in \{ -1, 1\}$ ,
\begin{equation}\label{3.1}
\begin{cases}\alpha_{\varepsilon} \in \mathbb{R}, \quad \alpha_{\varepsilon} \rightarrow \overline{\a}\in [ \overline{c}, 1 ] ,  \\
 a_{\varepsilon} \in \Omega, \quad \lambda_{\varepsilon} \in (0, \infty), \quad \lambda_{\varepsilon} d\left(a_{\varepsilon}, \partial \Omega\right) \rightarrow+\infty , \\
 v_{\varepsilon} \rightarrow 0  \text { in } H_{0}^{1}(\Omega), \quad v_{\varepsilon} \in E_{(a_{\varepsilon}, \lambda_{\varepsilon})}\end{cases}
\end{equation}
where $\overline{c}>0$ and for any $(a, \lambda) \in \Omega \times
(0, \infty) $, $ E_{(a, \lambda)}$ denotes the subspace of
$H_{0}^{1}(\Omega)$ defined by \eqref{vorthogonal}.
 \end{pro}
\begin{pf}
Let $u_\e$ be a solution of $(SC_\e)$ satisfying $(A)$ such that it
is weakly convergent to $0$ and does not converge strongly to $0$.
We denote
\begin{equation}\label{defw} w_\e := S_n^{n/4} u_\e / \| u_\e \| .
\end{equation}
We will prove that $(w_\e)$ is a Palais-Smale sequence of the functional $I_0$ through Proposition \ref{minim}. \\
From the definition of $w_\e$, it  follows that
$$ \| w_ \e \|^2  = S_n^{n/2} . $$
In addition, using Lemma \ref{lemma 3.1}, we have
$$ \| w_\e \|_{ L ^{2n / (n-2) }} = S _n^{ n/4} \frac{ \| u_\e \| _{  L ^{2n / (n-2) }}}{ \| u_\e \| } = S_n ^{ n/4} ( S_n^{-1/2} + o(1))  = S_n^{ (n-2) / 4 } + o(1), $$
which implies that
$$\displaystyle \lim_{\e \rightarrow 0}  \int_\O | w_\e |^{2n/(n-2)} = S_n^{n/2} . $$
Thus $ w_\e $ satisfies the assumptions of Proposition \ref{minim} and therefore $ w_\e$ is a Palais-Smale sequence for $ I_0$.
Furthermore,
$$ w_\e \rightharpoonup 0 \qquad \mbox{ and } \qquad  \| w_\e \|^2 = S_n ^{ n/2} . $$  Hence, we deduce from
\cite{Struwe} that
$$ w_\e = \gamma \d_{a_\e , \l _\e } + \wtilde{v}_\e  \qquad \mbox{ with } \quad    \gamma \in \{-1, 1\} , \quad \| \wtilde{v}_\e \| \to 0  \quad \mbox{ and } \quad  \l _\e d(a _\e , \partial \O ) \to \infty . $$
Following Bahri and Coron \cite{BC},  by modifying $ a_\e$ and $ \l_\e$, we can write $  w_\e $ as
\begin{equation}\label{decomposition}    \gamma  w_\e = \beta_\e P\d_{a_\e , \l_\e } + \wtilde{v} _ \e  \mbox{ with }  \beta_\e \to 1 , \, \, \l_\e d(a_\e, \partial \O ) \to \infty , \, \, \| \wtilde{v} _\e \| \to 0 \mbox{  and } \wtilde{v} _\e \in E_{ (a_\e , \l_\e)} .
\end{equation}
 Using  \eqref{defw},  \eqref{decomposition}, Lemma \ref{lemma 3.1}  and the fact that $ u_\e $ does not converge to $ 0 $, we get
$$ \g u_\e = \frac{ \| u_\e \|}{ S_n ^{ n/4 } }\g w_\e = \frac{ \| u_\e \|}{ S_n ^{ n/4 } } \beta_\e P\d_{a_\e , \l_\e } + \frac{ \| u_\e \|}{ S_n ^{ n/4 } } \wtilde{v} _ \e = \a_\e P\d_{a_\e , \l_\e } + v_\e $$
which completes  the proof.
\end{pf}


 \begin{rem}
Notice that $-u_\e$ is also a solution of $(SC_\e)$. Hence, in the sequel,  we will always assume that $u_\e$, solution of
$(SC_\e)$ satisfying the assumptions of the theorem, is written as
in \eqref{3.0}  with $ \gamma = 1 $  and \eqref{3.1} is satisfied.
\end{rem}

In the sequel, to simplify the notations, we set $\delta_{(a_{\varepsilon}, \lambda_{\varepsilon}) } = \delta_{\varepsilon} $ and $ P \delta_{a_{\varepsilon}, \lambda_{\varepsilon}} = P \delta_{\varepsilon} $.

Now we are going to study the concentration speed $\l_\e$ by giving
some relation between this parameter and $\e$ this will allow us to
expand $\ln (e+\d_\e)^\e$. This will be the subject of Lemmas
\ref{Lemma 4.2} and \ref{Lemma 3.3}. But some preliminaries results
are needed.


\begin{pro} \label{varphi} \cite{BCH, R} Let $a \in \Omega$ and $\lambda>0$ be such that $\lambda d(a, \partial \Omega)$ is large enough. For $\varphi_{(a, \lambda)}=\delta_{(a, \lambda)} - P \delta_{a, \lambda}$, we have the following estimates
 $$ 0 \leq \varphi_{(a, \lambda)} \leq \delta_{(a, \lambda)} , \leqno{(a)} $$
 $$ \varphi_{(a, \lambda)} = c_{0} \lambda^{\frac{2-n}{2}} H(a, .)+f_{(a, \lambda)}  , \leqno{(b)}  $$
where $f_{(a, \lambda)}$ satisfies
$$
f_{(a, \lambda)} = O\left(\frac{1}{\lambda^{ (n + 2 ) / 2}
d^{n}}\right) \quad \mbox{ and } \quad  \lambda \frac{\partial
f_{(a, \lambda)}}{\partial \lambda}=O\left(\frac{1}{\lambda^{ (n +2)
/ 2} d^{n}}\right),
$$
where $d$ is the distance $d(a, \partial \Omega)$,
$$  \left|\varphi_{(a, \lambda)}\right|_{L^{\frac{2
n}{n-2}}}  =O\left((\lambda d)^{\frac{2-n}{2}}\right),
\quad\left\|\varphi_{(a, \lambda)}\right\|=O\left((\lambda
d)^{\frac{2-n}{2}}\right), \quad \left|\lambda \frac{\partial
\varphi_{(a, \lambda)}}{\partial \lambda}\right|_{L^{\frac{2
n}{n-2}}} =O\left({(\lambda d)^{\frac{2-n}{2}}}\right). \leqno{(c)}
$$
\end{pro}
We also introduce the following.
\begin{lem}\label{glemmaB1}
Let $g_\e$ be the function defined on $ \R $ as follow
$$g_\e(U):=\ln(e + | U | )^\e.$$
\begin{enumerate}
\item For $\e$ small enough, and any $U\in \R$,
\begin{equation}\label{gB.1} |g'_\e(U)|\leq c\e \frac{1}{e +  | U | },\end{equation} and
\begin{equation}\label{gB''}
|g''_\e(U)|\leq c\e \frac{1}{(e + | U | )^2}.
\end{equation}
\item For any $\e>0$, and any $U\in \R$, we have
 \begin{equation}
|g_\e(U)-1|\leq \e \ln\ln(e + | U | ) \ln( e + | U |
)^\e.\end{equation}
\end{enumerate}
\end{lem}
\begin{pf} Notice that $ g_\e $ is an even function with respect to the variable $ U $.
\begin{enumerate}
\item For any $U\in (0, \infty)$, we have
\begin{equation}\label{g'}g'_\e(U)=\e\ln(e+U)^{\e-1}\frac{1}{e+U}
\end{equation} and
$$g''_\e(U)=\e(\e-1)\ln(e+U)^{\e-2}\frac{1}{(e+U)^2}-\e\ln(e+U)^{\e-1}\frac{1}{(e+U)^2}.$$
Thus, using the fact that $0<\e<1$, \eqref{gB.1} and \eqref{gB''}
follow.
\item A simple computation shows that we have
$$\frac{\partial (g_\e(U))}{\partial \e}= \ln\ln (e + | U | ) \, \ln(e + | U | )^\e $$
and by the mean value theorem we get the desired result.
\end{enumerate}
\end{pf}








It follows that

\begin{lem}\label{Lemma 4.2} Let $u_{\varepsilon}$ satisfy the assumption of Theorem \ref{1.4}. Then  $ \a_\e $ and $\lambda_{\varepsilon}$ occurring in \eqref{3.1} satisfy
$$
\a_\e^\frac{4}{n-2}\ln(\lambda_{\varepsilon})^{\varepsilon}
\rightarrow 1 \quad \text { as } \varepsilon \rightarrow 0.
$$
\end{lem}
\begin{pf}
Multiplying $(SC_\e)$ by $P\d_\e$ and integrating by parts over
$\O$, we get
\begin{equation}\label{2.14}
\a_\e \|P\d_\e\|^2 = \int_\O | u _\e | ^\frac{4}{n-2} u_\e \ln(e +
|u_\e | )^\e P\d_\e ,
\end{equation}
since $v_\e$ occurring in \eqref{3.1} is orthogonal to $P\d_\e$.\\
Recall that from \cite[Lemma 2.2]{BFG}, we have
\begin{equation}\label{4.19}
\|P\d_\e\|^2=S_n^{n /2}+o(1) \quad \text { as } \varepsilon
\rightarrow 0.
\end{equation}
Concerning the right hand side of \eqref{2.14}, let \be\label{O1O2}
\O_1:= \{ x \in \O : | v _\e | + | \a_\e \varphi_\e | \leq
\frac{1}{2} \a_\e  \d_\e \} \quad , \quad \O_2 := \O\setminus \O_1,
\ee where $ \varphi_\e := \varphi_{a_\e, \l_\e} $. Observe that, in
$ \O_2 $, it holds
$$ P\d_\e:= \d_\e - \varphi_\e \leq c ( | v_\e | + | \varphi_\e | ) \quad \mbox{ and } \quad  | u_\e | \leq c ( | v_\e | + | \varphi_\e | ) . $$
Therefore, using assumption $(A)$, Holder's inequality and Sobolev
embedding theorem, we obtain
\begin{align}\label{021} \int_{ \O_2} | u_\e | ^{ \frac{n+2}{n-2} } \ln ^\e ( e + | u_\e | ) P \d_\e  &
 \leq c \int_{ \O_2} | u_\e | ^{ \frac{4}{n-2} } \ln ^\e ( e + | u_\e | ) ( | v_\e |^2 + | \varphi_\e | ^2 ) \nonumber\\
 & \leq c  \Big( \int_{ \O} | u_\e | ^{ \frac{2n}{n-2} } \ln ^{\e n / 2} ( e + | u_\e | ) \Big)^{2/n} ( \| v _\e \|^2 + \| \varphi _\e \| ^2 )  = o(1) .
\end{align}
Regarding the integral over $ \O_1$, we write
\begin{align}\label{4.21}
 & \int_{\Omega_1} | u_{\varepsilon} |^{\frac{4}{n-2}} u_\e \ln(e + | u_\e | )^\e P\d_\e \nonumber\\
  &  =  \int_{\Omega_1} | u_{\varepsilon}|^{\frac{4}{n-2}} u_\e \Big[ g_\e(\alpha_{\varepsilon}  \delta_{\varepsilon} + (- \a_\e \varphi_\e+v_\e) ) - g_\e(\alpha_{\varepsilon}  \delta_{\varepsilon}) \Big] P \delta_{\varepsilon} + \int_{\Omega_1} | u_{\varepsilon} | ^{\frac{4}{n-2}} u_\e g_\e(\alpha_{\varepsilon}  \delta_{\varepsilon}) P \delta_{\varepsilon}
 : =I_1+I_2.
\end{align}
We claim that
\begin{equation}\label{4.22}
I_1=o(1).
\end{equation}
In fact, by the mean value theorem, there exists
$\theta=\theta(x)\in(0, 1)$ such that
$$ | I_1 | = \Big| \int_{\O_1}  | u_{\varepsilon}|^{\frac{4}{n-2}} u_\e g_\e'(\alpha_\e \d_\e + \theta ( - \a_\e \varphi_\e +  v_\e))  ( - \a_\e \varphi_\e + v_\e )   P\d_\e \Big| \leq c \e \int_{\O_1} | u_{\varepsilon}|^{\frac{n+2}{n-2}} \frac{  | - \a_\e  \varphi_\e + v_\e |   \d_\e }{ e + | \a_\e \d_\e + \theta ( - \a_\e \varphi _\e + v ) | } $$
where we have used \eqref{gB.1}. Observe that, in $ \O_1$, it holds
$$ \theta  | - \a_\e  \varphi_\e + v_\e |  \leq  | - \a_\e  \varphi_\e + v_\e | \leq \frac{1}{2} \a_\e \d_\e $$
which implies that
$$ | \alpha_\e \d_\e + \theta ( - \a_\e \varphi_\e +  v_\e) | \geq \alpha_\e \d_\e -  \theta | - \a_\e \varphi_\e +  v_\e | \geq \frac{1}{2} \a_\e \d_\e .  $$
Therefore, we get
$$ | I_1 |  \leq c\,  \e \int_{\O_1} | u_{\varepsilon}|^{\frac{n+2}{n-2}} {  | - \a_\e  \varphi_\e + v_\e |  }  \leq c\,  \e \int_{\O_1} | u_{\varepsilon}|^{\frac{n+2}{n-2}} ( |   \varphi_\e | +  | v_\e |  ) \leq c \, \e ( \| v_\e \| + \| \varphi_\e \| ) = o(1) .  $$
Hence our claim \eqref{4.22} is proved.

 Notice that
\begin{equation}\label{majln} 1 < \ln(e+\a_\e\d_\e)\leq c \ln\l_\e. \end{equation}
Recall that, we denoted by $ p:= (n+2) / (n-2)$. Using
\eqref{majln}, \eqref{3.1}, Holder inequality, Sobolev embedding
theorem and Proposition \ref{varphi}, we obtain
\begin{align}\label{4.231}
I_{2}
 & =\int_{\Omega_1} \left(\alpha_\e^{p}\d_\e^{p}+O\left(\d_\e^{p-1}(\varphi_\e+ | v_\e | ) + | v_\e | ^p\right)\right)\ln(e+\alpha_\e\d_\e)^\e (\d_\e-\varphi_\e)\nonumber\\
 & = \alpha_\e^{p}\int_{\Omega_1} \d_{\varepsilon}^{p+1}\ln(e + \alpha_\e\delta_{\varepsilon})^\e+O\left(\ln(\l_\e)^\e\int_\O\d_\e^p( | \varphi_\e | + | v_\e | )+\d_\e | v_\e | ^p\right)\nonumber\\
 & = \alpha_\e^{p}\int_{\Omega_1} \d_{\varepsilon}^{p+1}\ln(e + \alpha_\e\delta_{\varepsilon})^\e+O\left(\ln(\l_\e)^\e
 \left(\|\varphi_\e\|+\|v_\e\|+\|v_\e\|^p\right)\right)\nonumber\\
 & = \alpha_\e^{p}\int_{\Omega } \d_{\varepsilon}^{p+1}\ln(e + \alpha_\e\delta_{\varepsilon})^\e+o(\ln(\l_\e)^\e).
\end{align}
 Let $\Omega_\l:= B(a_\e, \l_\e^{-3/4})$. On one hand, using \eqref{majln}, we have
\begin{equation}\label{4.232}\int_{\O\setminus \O_\l}\delta_{\varepsilon}^{p+1}\ln (e+ \a_\e \delta_{\varepsilon})^\e\leq
  C \ln(\l _\e)^\e\int_{\O\setminus \O_\l}\delta_{\varepsilon}^{p+1}\leq \frac{C\ln(\l_\e)^\e}{\l_\e^{n/4}}=o(\ln(\l_\e)^\e).
  \end{equation}
  On the other hand, in $\O_\l$ we have $\l|x-a_\e|\leq \l_\e^{1/4}$ which implies that
  $$1+\l_\e^2|x-a_\e|^2\leq 1+\sqrt{\l_\e}\leq c\sqrt{\l_\e} \quad \hbox{ and } \quad c\l_\e^{(n-2)/4} \leq  \delta_{\varepsilon} \leq c_0\l_\e^{(n-2)/2} $$
  and therefore
  \be \label{X2} c \ln \l_\e \leq \ln (e + \a_\e \d_\e ) \leq c '  \ln \l_\e \quad \hbox{ and } \quad  \ln ^\e (\l_\e) (1 + O( \e) ) \leq  \ln ^\e (e + \a_\e \d_\e ) \leq  \ln ^\e( \l_\e ) (1 + O( \e) ) \quad \mbox{ in } \O_\l . \ee
  Therefore we deduce that
  \begin{align}\label{4.233}
  \int_{\O_\l} \delta_{\varepsilon}^{p+1}\ln (e+ \a_\e \delta_{\varepsilon})^\e & = \ln ^\e( \l_\e )  \int _{\O_\l} \delta_{\varepsilon}^{p+1} + O \Big( \e \ln ^\e( \l_\e )    \int_{\O_\l} \delta_{\varepsilon}^{p+1} \Big) \notag \\
   & = \ln ^\e( \l_\e )  (S_n^{n / 2}+o(1) ).
  \end{align}
 Combining \eqref{4.231}, \eqref{4.232} and \eqref{4.233}, we get
 \begin{equation}\label{4.23}
I_2= \a_\e^p\ln\left(\l_\e\right)^\e\left(S_n^{n /
2}+o(1)\right)+o(1).
\end{equation}
 Combining \eqref{2.14}, \eqref{4.19}, \eqref{021}-\eqref{4.22},
\eqref{4.23} and \eqref{3.1}, the lemma follows.
\end{pf}

Notice that, since $\a_\e$ occurring in \eqref{3.1} satisfies $\a_\e
\rightarrow \overline{\a}\in [\overline{c}, 1]$, we deduce that
\begin{equation}\label{3.2}
1\leq \ln (\l_\e)^\e\leq c.
\end{equation}
Next, as in \cite[Lemma 2.3]{BEGR}, we can easily prove the
following estimate :

\begin{lem} \label{Lemma 3.3} Equation \eqref{3.2} satisfied by the parameter $\l_\e$ implies that

\begin{enumerate}
\item  $\e\ln \ln \l_\e  \leq c $.
\item $\ln\left(e+\alpha_\varepsilon \delta_\varepsilon\right)^{\varepsilon}=
\ln(\lambda_\varepsilon^\frac{n-2}{2})^{\varepsilon }+ O\left(\e\ln
\left[\frac{\ln(e+\a_\e \delta_\e)}{\ln \l_\e
^\frac{n-2}{2}}\right]\right) \quad \text { in } \Omega $.
\item $\ln\left(e+\alpha_\varepsilon \delta_\varepsilon\right)^{\varepsilon}=
\ln(\lambda_\varepsilon^\frac{n-2}{2})^{\varepsilon }+\e\ln
\left[\frac{\ln(e+\a_\e \delta_\e)}{\ln \l_\e
^\frac{n-2}{2}}\right]+O\left(\e^2\ln ^2 \left[\frac{\ln(e+\a_\e
\delta_\e)}{\ln \l_\e ^\frac{n-2}{2}}\right] \right) \text { in }
\Omega.$
\end{enumerate}
\end{lem}
\begin{pf}
Claim $1$ is consequence of \eqref{3.2}.\\ Writing
$$ \ln^\e (e+ \a_\e \d_\e) = \ln^\e ( \l_\e ^{ (n-2)/2}) \mbox{exp} \Big( \e \ln \Big[ \frac{ \ln ( e+ \a_\e \d_\e) }{ \ln ( \l_\e ^{ (n-2)/2})} \Big] \Big) $$
and using \eqref{majln} and \eqref{3.2}, we deduce that \be
\label{X1} \e \, \Big| \ln \Big[ \frac{ \ln ( e+ \a_\e \d_\e) }{ \ln
( \l_\e ^{ (n-2)/2})} \Big] \Big| \leq c \qquad \mbox{ uniformly in
} \O . \ee Thus, by Taylor expansion of $ e^t$ and using
\eqref{3.2}, Claims $2$ and $3$ follow.
\end{pf}

We are now able to study the $v_{\varepsilon}$-part of
$u_{\varepsilon}$.

\section{Estimating $v_\e$}
\begin{lem}\label{integralv}Let $u_{\varepsilon}$ satisfy the assumption of Theorem \ref{1.4}. Then $v_{\varepsilon}$ occurring in \eqref{3.1} satisfies

$$  \int_\O | v_\e |^{2n/(n-2) } \ln ^{\e n /2}  ( e + | v_\e | )  = o(1) . $$
\end{lem}
\begin{pf}
To estimate this integral, we need the following result.
\begin{lem}\label{lem1}
Let $$ \psi_\e(t) := | t|^{2n / (n-2) } \ln ^{\e n /2} (e + | t | )
\qquad t \in \R. $$ For each $ t , s \in \R$, it holds
$$ | \psi_\e ( t+ s ) - \psi_\e (t) - \psi_\e (s) | \leq c | t |^{(n+2) / (n-2) } \ln^{\e n / 2 } ( e +  | t | ) | s | + c | s |^{(n+2) / (n-2) } \ln^{\e n / 2 } ( e +  | s | ) | t | . $$
\end{lem}
Applying Lemma \ref{lem1} and using the fact that $P\d_\e \leq
\d_\e$, we get
\begin{align}\label{v0}\int_\O \psi_\e( \a P\d_\e + v_\e ) =& \int_\O \psi_\e( \a_\e P\d_\e  )  + \int_\O \psi_\e(  v_\e )  \nonumber\\&+ O \Big( \int_\O \d_\e^{(n+2)/(n-2) } \ln ^{\e n / 2} ( e +  \d_\e ) | v_\e | +  \int_\O | v_\e |^{(n+2)/(n-2) } \ln ^{\e n / 2 } ( e +  | v_\e | ) \d_\e \Big) .
\end{align}
By assumption $(A)$, the left integral is
\begin{align}\label{v1}
\int_\O \psi_\e( \a_\e P\d_\e + v_\e )=S_n ^{n/2} (1+o(1)).
\end{align}
We are going to estimate each term on the right hand-side in the above equality.\\
The first term will be computed as the integral $I_2$ in
\eqref{4.21}. By the mean value theorem and taking into account
\eqref{3.1}, \eqref{gB.1}, \eqref{majln}, \eqref{3.2} and the fact
that
\begin{equation}\label{pdless}
P\d_\e \leq \d_\e - t  \varphi_\e, \quad\hbox{for each } t \in [0,1],
 \end{equation}
 we have
\begin{align}\label{4.234}
\int_\O \psi_\e( \a_\e P\d_\e  )=&
 \int_{\Omega} \a_{\varepsilon}^{p+1}P\d_\e^{p+1}g_\frac{\e n}{2}(\a_\e\d_\e)+
 \int_{\Omega} \a_{\varepsilon}^{p+1}P\d_\e^{p+1}\left[g_\frac{\e n}{2}(\a_\e P\d_\e)-g_\frac{\e n}{2}(\a_\e\d_\e)\right]
\nonumber\\
=&\int_{\Omega}
\left(\alpha_\e^{p+1}\d_\e^{p+1}+O\left(\d_\e^{p}\varphi_\e\right)\right)\ln(e+\alpha_\e\d_\e)^{\frac{\e
n}{2}}
+O\left(\e\int_{\Omega} \d_{\varepsilon}^{p}\frac{\a_\e P\d_\e}{e+\alpha_\e\d_\e-\theta \a_\e\varphi_\e}\varphi_\e\right)\nonumber\\
=&\alpha_\e^{p+1}\int_{\Omega} \d_{\varepsilon}^{p+1}\ln(e + \alpha_\e\delta_{\varepsilon})^{\frac{\e n}{2}}+O\left(\e+\int_\O\d_\e^p\varphi_\e\right) \nonumber\\
=&\alpha_\e^{p+1}\int_{\Omega} \d_{\varepsilon}^{p+1}\ln(e + \alpha_\e\delta_{\varepsilon})^{\frac{\e n}{2}}+O\left(\e+\frac{1}{(\l_\e d_\e)^\frac{n-2}{2}}\right)\nonumber\\
=&\alpha_\e^{p+1}\int_{\Omega} \d_{\varepsilon}^{p+1}\ln(e +
\alpha_\e\delta_{\varepsilon})^{\frac{\e n}{2}}+o(1).
\end{align}
Proceeding as in \eqref{4.232} and \eqref{4.233}, by splitting the
integral over $\O_\l$ and $\O_\l^c$, and substituting $\ln(e+ \a_\e
\d_\e)^\e$ for $\ln(e+ \a_\e\d_\e)^{\e n /2}$, we obtain in the same
way
\begin{equation}\label{4.235}\int_{\Omega} \d_{\varepsilon}^{p+1}\ln(e + \alpha_\e\delta_{\varepsilon})^{\frac{\e n}{2}}=\ln\left(\l_\e\right)^\frac{\e n}{2}\left(S_n^{n / 2}+o(1)\right).\end{equation}
Combining \eqref{4.234} and \eqref{4.235}, together with the result
of Lemma \ref{Lemma 4.2}, we derive that
\begin{equation}\label{v2}
\int_\O \psi_\e( \a_\e P\d_\e  )=S_n^{n / 2}+o(1).
\end{equation}
 Now,  using \eqref{majln} and \eqref{3.2}, observe  that
\begin{equation}\label{v3} \int_\O \d_\e^{(n+2)/(n-2) } \ln ^{\e n / 2} ( e +  \d_\e ) | v_\e | \leq c \int_\O \d_\e^{(n+2)/(n-2) } | v_\e | \leq c \| v_\e \| = o(1)  .\end{equation}
It remains to estimate the last term in the right hand side. We will
split it into two integrals. \\For $ M$ a fixed large positive
constant, we denote
$$ \O_1 := \{ x \in \O : | v_\e | \leq M \d_\e \} \qquad \mbox{ and } \qquad \O_2 := \O \setminus \O_1 .$$
It follows that, in $ \O_2$, the function $ \d_\e$ is small with
respect to $ | v_\e |$. However,  in $ \O_1$ through \eqref{majln},
\eqref{3.2}, we have
$$ \ln^{\e n / 2 }  (e + |v_\e|) \leq \ln^{\e n / 2 }  (e + M \d_\e) \leq \ln^{\e n / 2 }  ( M(e+ \d_\e )) \leq ( c \ln (\l_\e^{ (n-2) / 2 }))^{ \e n / 2 }  \leq c . $$
Thus, we obtain
\begin{align}\label{v4}
\int_\O | v _\e|^{(n+2)/(n-2) } \ln ^{\e n / 2 } ( e +  | v_\e | ) \d_\e & \leq \frac{1}{M} \int_{\O_2} | v_\e |^{ 2 n /(n-2) } \ln ^{\e n / 2 } ( e +  | v_\e | ) + c \int_{\O_1} | v_\e |^{(n+2)/(n-2) }  \d_\e \nonumber\\
 & \leq \frac{1}{M} \int_\O \psi_\e(  v_\e )   + o(1)
 \end{align}
where we have used Holder's inequality, Sobolev
embedding theorem and the fact that $\|v_\e\|=o(1)$ as $\e$ goes to zero.\\
Hence \eqref{v0}, \eqref{v1}, \eqref{v2}, \eqref{v3} and \eqref{v4}
lead to
$$  \Big( 1 + O \Big( \frac{1}{M} \Big) \Big) \int_\O \psi_\e(  v_\e )  = o(1) $$
which concludes the proof of this lemma.
\end{pf}


\begin{lem}\label{Lemma 4.6} Let $u_{\varepsilon}$ satisfy the assumption of Theorem \ref{1.4}. Then $v_{\varepsilon}$ occurring in \eqref{3.1} satisfies
\begin{align}
\left\|v_{\varepsilon}\right\| \leq & C
\frac{\varepsilon}{\ln\lambda_\varepsilon}+C\left\{\begin{array}{ll}
\frac{1}{(\l_\e d_\e)^{n-2}} & \hbox{if }n<6,  \\
                                                      \frac{\ln(\l_\e d_\e)}{(\l_\e d_\e)^4} & \hbox{if }n=6,  \\
                                                     \frac{1}{(\l_\e d_\e)^{(n+2)/2}} & \hbox{if }n>6,
                                                     \end{array}\right.
\end{align}
with $C$ independent of $\varepsilon$.
\end{lem}
\begin{pf}  Multiplying $\left(SC_\e\right)$ by $v_{\varepsilon}$ and integrating on $\Omega$, we obtain
\begin{equation}\label{vest0}
\int_{\Omega} \nabla u_{\varepsilon} \cdot \nabla v_{\varepsilon} -
\int_{\Omega}  | u_{\varepsilon}| ^{4/(n-2) } u_\e \ln(e+ | u_\e |
)^\e v_\e = 0 .
\end{equation}
Let \be \label{Ov1} \O_v := \{ x \in \O : \a_\e  P \d _\e \leq  2| v
_\e | \} . \ee In $ \O_v $, it holds that \be \label{X3} e + | u_\e
| \leq e + 3 | v_\e | \leq 3 ( e + | v_\e| ) \quad  \hbox{ and }
\quad  \ln (   e + | u_\e | ) \leq \ln 3 + \ln ( e + | v_\e| ) \leq
c \ln ( e + | v_\e| ) . \ee Thus, we get $$   \ln^\e  (   e + | u_\e
| ) \leq c  \ln^\e  (   e + | v_\e | )  .$$ Therefore, using
Holder's inequality, Sobolev embedding theorem and Lemma
\ref{integralv}, we obtain
\begin{align}
\label{vest1}\int_{ \O_v} | u_\e |^{ p}  \ln ^\e ( e+ | u_\e | )  |
v_\e | &\leq c  \int_{ \O_v} | v_\e |^{ p+1}  \ln ^\e ( e+ | v_\e |
) \leq c   \int_{ \O }| v_\e |^2 | v_\e |^{ p-1}  \ln ^\e ( e+ |
v_\e | )
\nonumber\\
&\leq c \Big( \int_{ \O }| v_\e |^{2n/(n-2)} \Big)^{(n-2)/n}  \Big(
\int_{ \O }| v_\e |^{2n/(n-2)} \ln ^{ \e n / 2 } ( e + | v_\e |
)\Big)^{ 2 / n } = o( \| v_\e \| ^2 ) .
\end{align}
 Using \eqref{3.1}, \eqref{vest0}, \eqref{vest1} and the function $g_\e$ introduced in Lemma \ref{glemmaB1}, we have
\begin{align}\label{estv0}
&\int_{\Omega}\left|\nabla v_{\varepsilon}\right|^{2}-
\int_{\Omega\setminus \O_v} | u_\e | ^{\frac{ 4}{n-2} } u_\e
\left[g_\e(u_\e)-g_\e(\alpha_\e P\delta_\e)\right] v_{\varepsilon} -
\int_{ \O \setminus \Omega_v }  | u_\e | ^{ \frac{4}{n-2}} u_\e
g_\e(\alpha_\e P\delta_\e)v_\e +o( \| v_\e \| ^2 ) =0.
\end{align}
Observe that, for each $ t \in [0,1]$, it holds \be
\label{vlesspdelta} | v_\e | \leq  \frac{1}{2}  \a_\e P \d_\e  \leq
\a_\e P \d_\e - t | v_\e | \leq  \a_\e P \d_\e + t v_\e \leq
\frac{3}{2} \a_\e P \d_\e  \quad \mbox{  in } \O\setminus \O_v. \ee
By the mean value theorem, Holder's inequality, Sobolev embedding
theorem and Lemma \ref{lemma 3.1} and using \eqref{gB.1} and
\eqref{vlesspdelta}, we have
\begin{align}\label{estv1}
\int_{\O\setminus \O_v} | u_\e | ^{\frac{ 4}{n-2} } u_\e\left[g_\e(u_\e)-g_\e(\alpha_\e P\delta_\e)\right] v_{\varepsilon}=& \int_{\O\setminus \O_v}
 | u_\e | ^{\frac{ 4}{n-2} } u_\e g'_\e(\alpha_\e P\delta_\e+\theta v_\e) v_\e^2 \nonumber\\
=&O\left(\e\int_{\O\setminus \O_v}u_\e^{p-1}\frac{\alpha_\e P\delta_\e + v_\e }{e +\alpha_\e P\delta_\e+\theta v_\e}v_\e^2  \right)\nonumber\\
=&O\left(\e \int_{\Omega} | u_\e | ^{p-1}v_\e^2\right)\nonumber\\
=&O\left(\e \|v_\e\|^2\right).
\end{align}
Moreover, using the fact that $ P \d_\e \leq \d_\e $, \eqref{majln},
\eqref{3.2} and following the proof of \eqref{vest1}, we get
\begin{align}\label{estv3}
 & \int_{\Omega \setminus \O_v } | u_\e | ^{p-1 } u_\e g_\e(\alpha_\e P\delta_\e) v_\e \notag \\
& =  \a_\e^p\int_{\Omega \setminus \O_v  }P\d_\e^p g_\e(\alpha_\e P\delta_\e) v_\e+p\a_\e^{p-1}\int_{\Omega \setminus \O_v  }P\d_\e^{p-1} g_\e(\alpha_\e P\delta_\e) v_\e^2+O\left(\int_{\Omega \setminus \O_v  } \d_\e^{p-2} | v_\e | ^3 + \int_{\Omega \setminus \O_v  }  | v_\e | ^{p+1}\right)\nonumber\\
 & = \a_\e^p\int_{\Omega}P\d_\e^p g_\e(\alpha_\e P\delta_\e)
v_\e+p\a_\e^{p-1}\int_{\Omega}P\d_\e^{p-1} g_\e(\alpha_\e
P\delta_\e) v_\e^2+o\left(\| v_\e\|^2\right).
\end{align}
Let
$$ Q_\e(v, v):=\|v\|^{2}-p \int_{\Omega}(\a_\e P \d_\e )^{p-1} g_\e(\alpha_\e P\delta_\e)v^2 \quad \hbox{ and } \quad L_\e(v):=\int_{\Omega}(\a_\e P \d_\e )^p g_\e(\alpha_\e
P\delta_\e)v .$$ From \eqref{estv0}, \eqref{estv1} and
\eqref{estv3}, we get
\begin{equation}\label{vest0'}
Q\left(v_{\varepsilon},
v_{\varepsilon}\right)-L_{\varepsilon}\left(v_{\varepsilon}\right)+R_\e(v_\e)=0
\end{equation}
where $ R_\e(v)$ is a $\mathcal{C}^2$ function satisfying
\begin{equation}R_{\eps}(v_\e)=o(\|v_\e\|^{2})
\end{equation}
uniformly with respect to $\e,\a_\e,\l_\e,a_\e$ for $\e$ small enough and $\a_\e,\l_\e,a_\e$ verify \eqref{3.1}.\\
By the mean value theorem and using \eqref{gB.1}, \eqref{majln},
\eqref{3.2} and Lemma  \ref{Lemma 3.3}, we find
\begin{align}
Q_\e(v_\e,v_\e) & = \|v\|^{2}-p \int_{\Omega}(\a_\e  P\d_\e )^{p-1} g_\e(\alpha_\e \delta_\e)v^2 + p \int_{\Omega}(\a_\e  P\d_\e )^{p-1} g_\e'(\alpha_\e \d_\e - \theta \alpha_\e  \varphi_\e)\a_\e\varphi_\e v_\e ^2\nonumber\\
 & = \|v_\e\|^2 - \a_\e^{p-1}\ln(\l^\frac{n-2}{2})^\e p\int_{\Omega } \d_\e^{p-1} v_\e^2 + O\left( \e\int_\O \d_\e^{p-1}\ln
\left[\frac{\ln(e+\a_\e \delta_\e)}{\ln \l_\e ^\frac{n-2}{2}}\right]
v_\e^2+\int_\O \d_\e^{p-2}\varphi_\e v_\e^2\right)\nonumber\\
& \quad + O\left(\e\int_{\Omega}\d_\e^{p-1}\frac{\alpha_\e P\d_\e
}{e +  | \alpha_\e \d_\e -  \theta \alpha_\e  \varphi_\e |
}v_\e^2\right) .  \label{3.18bis}
\end{align}
 Notice that, using Proposition \ref{varphi}, we get
$$ \int_\O \d_\e^{p-2}\varphi_\e v_\e^2 \leq c \int_\O \d_\e^{3/(n-2)} \varphi_\e ^{1/(n-2)}  v_\e^2 \leq c \| \varphi_\e \|^{1/(n-2)} \| v _\e \| ^2 = o( \| v _\e \|^2) . $$
In addition,  taking into account \eqref{pdless},
we get
$$ \e\int_{\Omega}\d_\e^{p-1}\frac{\alpha_\e P\d_\e  }{e +  | \alpha_\e \d_\e -  \theta \alpha_\e  \varphi_\e | }v_\e^2 \leq \e\int_{\Omega}\d_\e^{p-1} v_\e^2 \leq c \, \e \, \| v _\e\| ^2  = o(  \| v _\e\| ^2 ) .  $$
Concerning the other integral in the remainder terms of
\eqref{3.18bis}, as in the proof of \eqref{4.231},  let $ \O_\l :=
B(a_\e , \l_\e ^{-3/4}) $, using \eqref{X1} and the first claim of
\eqref{X2}, we have
\begin{align*}
 & \e \, \Big| \int_{\O \setminus \O_\l} \d_\e^{p-1}\ln \left[\frac{\ln(e+\a_\e \delta_\e)}{\ln \l_\e ^\frac{n-2}{2}} \right] v_\e^2 \Big| \leq c \int_{\O \setminus \O_\l} \d_\e^{p-1} v_\e^2 \leq c \| v _\e \|^2 \Big( \int_{\O \setminus \O_\l} \d_\e^{p + 1} \Big)^{ 2/n} = o( \| v_\e\|^2 ) , \\
 &  \e \, \Big| \int_{ \O_\l} \d_\e^{p-1}\ln \left[\frac{\ln(e+\a_\e \delta_\e)}{\ln \l_\e ^\frac{n-2}{2}} \right] v_\e^2 \Big| \leq c  \e \, \Big| \int_{ \O_\l} \d_\e^{p-1} v_\e^2 \Big| \leq c \, \e \| v_\e \|^2 \Big( \int_{\O } \d_\e^{p + 1} \Big)^{ 2/n} \leq  c \, \e \| v_\e \|^2 = o  ( \| v_\e \|^2 ) .
  \end{align*}
  Hence, using Lemma \ref{Lemma 4.2}, \eqref{3.18bis} becomes
\be \label{3.18} Q_\e(v_\e,v_\e) =  Q_0(v_\e,v_\e)+ o(\|v_\e\|^2)
\quad \mbox{ where } \quad  Q_ 0( v,v):= \|v\|^2- p\int_{\Omega
}\d_\e^{p-1} v^2 . \ee
 According to \cite{B1,R}, $Q_0$ is coercive
in the space  $ E_{(a_{\varepsilon}, \lambda_{\varepsilon}) }$, that
is, there exists some constant $c>0$ independent of $\varepsilon$,
for $\varepsilon$ small enough, such that

\begin{equation}\label{vest1'}
Q_0(v, v) \geq c\|v\|^{2} \quad \forall v \in
E_{\left(a_{\varepsilon}, \lambda_{\varepsilon}\right)}.
\end{equation}

In the sequel, we study $L_{\e}(v_\e)$. We claim that

\begin{align}\label{vest2'}
L_{\e}(v_\e) =&\int_{\Omega}\a_\e^p\d_\e^p
g_\e\left(\a_\e\delta_{\varepsilon}\right)v_\e
+\left\{\begin{array}{ll}
O\left(\frac{\|v\|}{(\l_\e d_\e)^{n-2}}\right) & \hbox{if }n<6 ,  \\
                                                      O\left(\|v\|\frac{\ln(\l_\e d_\e)}{(\l_\e d_\e)^4}\right) & \hbox{if }n=6 , \\
                                                      O\left(\frac{\|v\|}{(\l_\e d_\e)^{(n+2)/2}}\right) & \hbox{if }n>6.
                                                     \end{array}\right.
\end{align}
In fact, by the mean value theorem, Holder's inequality  and using
\eqref{gB.1}, \eqref{majln}, \eqref{3.2}, \eqref{pdless} and Lemma
\ref{varphi}, we get
\begin{align}\label{1vest2'}
L_{\e}(v_\e)
 & =  \int_{\Omega}\a_\e^pP\d_\e^p g_\e\left(\a_\e\delta_{\varepsilon}\right)v_\e-\int_{\Omega}\a_\e^pP\d_\e^pg_\e ' \left(\alpha_\e \d_\e - \alpha_\e \theta \varphi_\e\right)\a_\e\varphi_\e v_\e \nonumber\\
 & =  \int_{\Omega}\a_\e^p\d_\e^p g_\e\left(\a_\e\delta_{\varepsilon}\right)v_\e+O\left(\int_{\Omega} \delta_{\varepsilon}^{p-1} \varphi_{\varepsilon}|v_\e|+\e\int_{\Omega} \delta_{\varepsilon}^{p-1} \varphi_{\varepsilon}|v_\e|\right)\nonumber\\
 & = \int_{\Omega}\a_\e^p\d_\e^p g_\e\left(\a_\e\delta_{\varepsilon}\right)v_\e +O\left(\int_{B_\e}\d_\e^{p-1}\varphi_\e | v_\e | + \int_{B_\e ^c }\d_\e^{p} | v_\e | \right)\nonumber\\
& = \int_{\Omega}\a_\e^p\d_\e^p
g_\e\left(\a_\e\delta_{\varepsilon}\right)v_\e
+O\left(\|v_\e\|\left(\frac{1}{\l_\e^{(n-2)/2}d_\e^{n-2}} \Big(
\int_{B_\e}\d_\e^\frac{8n}{n^2-4}
\Big)^\frac{n+2}{2n}+\frac{1}{(\l_\e d_\e)^{(n+2)/2}}\right)\right),
\end{align}
where $B_\e$ denotes the ball of center $a_\e$ and radius $d_\e/2$.
Furthermore, we have the following computation
\begin{equation}\label{2vest2'}
\Big( \int_{B_\e}\d_\e^\frac{8n}{n^2-4} \Big)^\frac{n+2}{2n}
=\left\{\begin{array}{ll}
O\left(\frac{1}{(\l_\e d_\e)^{(n-2)/2}}\right) & \hbox{if }n<6, \\
                                                      O\left(\frac{\ln^{2/3}(\l_\e d_\e)}{\l_\e ^2}\right) & \hbox{if }n=6, \\
                                                      O\left(\frac{d_\e^{(n-6)/2}}{\l_\e^2}\right) & \hbox{if }n>6.
                                                     \end{array}\right.
\end{equation}
Hence, combining \eqref{1vest2'} and \eqref{2vest2'}, Claim \eqref{vest2'} follows.\\
Lastly, we estimate the integral defined in \eqref{vest2'}. By using
Lemma \ref{Lemma 3.3}, \eqref{3.2}  and
the fact that $v_\e \in E_{\left(a_{\varepsilon},\l_\e\right)}$, we
get
\begin{align}\label{1v}
\int_{\Omega}\a_\e^p\d_\e^p
g_\e\left(\a_\e\delta_{\varepsilon}\right)v_\e=& \a_\e^p\ln
(\l_\e^\frac{n-2}{2})^\e\int_{\Omega}\delta_{\varepsilon}^pv_\e+
O\left(\e\int_{\O}\d_\e^p|v_\e| \Big| \ln \left[\frac{\ln(e+\a_\e \delta_\e)}{\ln \l_\e ^\frac{n-2}{2}}\right] \Big| \right)\nonumber\\
=&O\left(\e\int_{\O}\d_\e^p|v_\e| \Big| \ln \left[\frac{\ln(e+\a_\e
\delta_\e)}{\ln \l_\e ^\frac{n-2}{2}}\right] \Big| \right) .
\end{align}
Note that \eqref{majln} implies
 \begin{equation}\label{lnln}\left| \ln \left[\frac{\ln(e+\a_\e \delta_\e(x))}{\ln \l_\e ^\frac{n-2}{2}}\right]\right|\leq C \ln\ln \l_\e, \qquad  \forall x\in
 \Omega.
 \end{equation}

As in the proof of \eqref{4.231}, let $\O_\l:=
B(a_\e,\l_\e^{-3/4})$. We split the integral over $\O_\l$ and $\O
\setminus \O_\l$. On one hand, using \eqref{lnln}, Holder inequality
and Sobolev embedding theorem, we have
$$
 \int_{\O\setminus \O_\l}\d_\e^p |v_\e| \left|  \ln \left[\frac{\ln(e+\a_\e \delta_\e)}{\ln \l_\e
^\frac{n-2}{2}}\right] \right|  \leq C \ln\ln \l_\e
\left(\int_{\O\setminus \O_\l}\d_\e^{p+1}\right)^\frac{n+2}{2n}
\|v_\e\|
 \leq  C\frac{\ln \ln \l_\e}{ \l_\e^\frac{n+2}{8}} \|v_\e\|.
$$
 Thus \begin{equation}\label{3V}\e\displaystyle \int_{\O\setminus \O_\l}\d_\e^p |v_\e| \left|
 \ln \left[\frac{\ln(e+\a_\e \delta_\e)}{\ln \l_\e ^\frac{n-2}{2}}\right] \right|  =  O\left(\frac{\e}{\ln \l_\e}\|v_\e\|\right).\end{equation}
 On the other hand, in $\O_\l $ we have $\l_\e^2|x-a_\e|^2\leq \l_\e^{1/2}$ and
  $\frac{\l_\e}{1+\l_\e^2|x-a_\e|^2}\geq c\sqrt{\l_\e}$ which implies
 \be \label{a2} \d_\e(x)\geq c \l_\e^ {(n-2) /4} \quad \forall \, x\in \O_\l. \ee
  Since $e+\a_\e\d_\e \geq c \l_\e^ {(n-2)/ 4} $ we get $\ln (e+\a_\e\d_\e) \geq \ln(c)+\frac{1}{2} \ln \l_\e^{(n-2)/2} $ and therefore
  $$ \frac{\ln (e+\a_\e\d_\e)}{\ln \l_\e^ {(n-2) /2}}\geq \frac{1}{2}+ \frac{\ln(c)}{\ln\l_\e^ {(n-2) / 2}}\geq \frac{1}{4}. $$
  Now, since $ e+\a_\e\d_\e \leq c \l_\e^ { (n-2) / 2} $ we get $\ln (e+\a_\e\d_\e) \leq \ln(c)+ \ln \l_\e^ {(n-2) / 2} $ and therefore
  $$\frac{\ln (e+\a_\e\d_\e)}{\ln \l_\e^{(n-2)/2}}\leq 1+ \frac{\ln(c)}{\ln\l_\e^{(n-2)/2}}\leq \frac{3}{2}. $$
  Using the fact that $|\ln (t)|\leq c |t-1|$ for all $t\in [1/4,3/2]$, we get
  \begin{align}\label{1ln}\left|\ln \left[\frac{\ln(e+\a_\e \delta_\e)}{\ln \l_\e ^ {(n-2)/2}}\right]\right|\leq & c \left|\frac{\ln(e+\a_\e \delta_\e)}{\ln \l_\e ^ { (n-2) / 2}}
  -1\right|\nonumber\\
  \leq & \frac{c}{\ln\l_\e^\frac{n-2}{2}}\left|\ln \left(\frac{e}{\l_\e^ {(n-2)/2}}+
  \frac{\a_\e c_0}{(1+\l_\e^2|x-a_\e|^2)^ {(n-2)/2}}\right)\right|, \, \forall x\in B_\l.
  \end{align}
  Using \eqref{1ln}, Holder inequality and Sobolev
embedding theorem, we have
  \begin{align}\label{4V}
  \e\displaystyle \int_{ B_\l}\d_\e^p |v_\e|\ln \left[\frac{\ln(e+\a_\e \delta_\e)}{\ln \l_\e ^\frac{n-2}{2}}\right]\leq&
c \frac{\e}{\ln \l_\e}\|v_\e\|\left(\int_{B_\l}\d_\e^{p+1}\left|\ln
\left(\frac{e}{\l_\e^\frac{n-2}{2}}+
 \frac{\a_\e c_0}{(1+\l_\e^2|x-a_\e|^2)^\frac{n-2}{2}}\right)\right|^\frac{2n}{n+2}\, dx\right)^\frac{n+2}{2n}\nonumber\\
\leq& c \frac{\e}{\ln
\l_\e}\|v_\e\|\left(\int_{\tilde{B}_\l}\frac{1}{(1+|y|^2)^n}\left|\ln
\left(\frac{e}{\l_\e^\frac{n-2}{2}}+ \frac{\a_\e
c_0}{(1+|y|^2)^\frac{n-2}{2}}\right)\right|^\frac{2n}{n+2}\,
dy\right)^\frac{n+2}{2n}, 
  \end{align}
  where $\tilde{B}_\l=B(0, \l_\e^ {1/4})$ and we have used the change of coordinates
  $y=\l_\e(x-a_\e)$. We claim that
  \begin{equation}\label{5V}
  K:=\int_{\tilde{B}_\l}\frac{1}{(1+|y|^2)^n}\left|\ln
\left(\frac{e}{\l_\e^{(n-2)/2}}+ \frac{\a_\e
c_0}{(1+|y|^2)^{(n-2)/2}}\right)\right|^\frac{2n}{n+2}\, dy \leq C.
  \end{equation}
 Indeed,  observe that, in $\tilde{B}_\l$ we have $1+ |y|^2\leq 2\sqrt{\l_\e}$  and thus $\frac{\a_\e c_0}{(1+|y|^2)^\frac{n-2}{2}}\geq \frac{c}{\l_\e^\frac{n-2}{4}}
  >>\frac{e}{\l_\e^\frac{n-2}{2}}$. So, we get
  $$\frac{\a_\e c_0}{(1+|y|^2)^ {(n-2)/2}}\leq\frac{e}{\l_\e ^ {(n-2)/2}} +\frac{\a_\e c_0}{(1+|y|^2)^ {(n-2)/2}}\leq \frac{C}{(1+|y|^2)^ {(n-2)/2}}.$$
  Hence $$0\leq \ln \left(\frac{e}{\l_\e ^ {(n-2)/ 2}} +\frac{\a_\e c_0}{(1+|y|^2)^ {(n-2 )/ 2}}\right)-
  \ln \left(\frac{\a_\e c_0}{(1+|y|^2)^ {(n-2)/2}}\right)\leq C,\, \forall y \in \tilde{B}_\l.$$
The last inequality implies that
\begin{align*}
K\leq &C \int_{\tilde{B_\l}}\frac{1}{(1+|y|^2)^n}+\frac{1}{(1+|y|^2)^n}\ln \left(\frac{\a_\e c_0}{(1+|y|^2)^\frac{n-2}{2}}\right)^\frac{2n}{n+2}\, dy\nonumber\\
\leq& C\int_{\mathbb{R}^n}\frac{1}{(1+|y|^2)^n}+\frac{1}{(1+|y|^2)^n}\ln \left(\frac{\a_\e c_0}{(1+|y|^2)^\frac{n-2}{2}}\right)^\frac{2n}{n+2}\, dy\nonumber\\
\leq& C
\end{align*}
and the claim follows. Through \eqref{3V}, \eqref{4V} and
\eqref{5V}, we obtain
\begin{equation}\label{2V}
\e\int_{\O}\d_\e^p|v_\e|\ln \left[\frac{\ln(e+\a_\e \delta_\e)}{\ln
\l_\e
^\frac{n-2}{2}}\right]=O\left(\|v_\e\|\frac{\e}{\ln\l_\e}\right).
\end{equation}
\eqref{1v} and \eqref{2V} imply that
\begin{equation}\label{6V}
\int_{\Omega}\a_\e^p\d_\e^p
g_\e\left(\a_\e\delta_{\varepsilon}\right)v_\e
=O\left(\|v_\e\|\frac{\e}{\ln\l_\e}\right).
\end{equation}
Combining \eqref{vest0'}-\eqref{vest2'} and \eqref{6V}, we obtain
the desired estimate.
\end{pf}

\section{Proof of Theorem \ref{1.4}}

We start this section by the following lemma which will be useful
later.

 \begin{lem} \label{a34}Let $ d_\e := d(a_\e , \partial \O)$,  $ \eta_\e := \min( d_\e , \l_\e ^{ -3/4} )$ and $ B_\eta := B( a_\e , \eta_\e)$. For $\e$ small enough, we have
 \begin{align}
 &  \frac{\ln(e+\a_\e \delta_\e)}{\ln \l_\e ^\frac{n-2}{2}}
  -1=\frac{\ln(\a_\e c_0)}{\ln \l_\e ^\frac{n-2}{2}}-\frac{n-2}{2}\frac{\ln(1+\l_\e^2|x-a_\e|^2)}{\ln \l_\e ^\frac{n-2}{2}}
  +O\left(\frac{1}{\l_\e^\frac{n-2}{4}\ln \l_\e
  ^\frac{n-2}{2}}\right) \qquad \mbox{ in } B_\eta ,\label{a3} \\
 &  \frac{\ln(e+\a_\e \delta_\e)}{\ln \l_\e ^\frac{n-2}{2}}
  -1=O\left(\frac{1+\ln(1+\l_\e^2|x-a_\e|^2)}{\ln \l_\e
  ^\frac{n-2}{2}}\right) \qquad  \hbox{ in } B_\eta. \label{a4}
 \end{align}
 \end{lem}
 \begin{pf}
 Since $ B_\eta$ is a subset of $ \O _\l $ (which is introduced in the proof of Lemma \ref{Lemma 4.6}) we derive that \eqref{a2} holds true in $ B_\eta $. In addition, for $ x \in  B_\eta$ we have
 \begin{align*}
 \ln(e+\a_\e \delta_\e(x)) = & \ln \left(\a_\e
 \d_\e (x) \left[1+\frac{e}{\a_\e\d_\e(x)}\right]\right)\\
  =&\ln(\a_\e c_0)+\ln \l_\e
  ^\frac{n-2}{2}-\frac{n-2}{2}\ln(1+\l_\e^2|x-a_\e|^2)+\ln\left(1+\frac{e}{\a_\e\d_\e(x)}\right).
 \end{align*}
 Therefore, using \eqref{a2}, we get
 \begin{align*}
 \frac{\ln(e+\a_\e \delta_\e (x) )}{\ln \l_\e ^\frac{n-2}{2}}
  -1=&\frac{\ln(\a_\e c_0)}{\ln \l_\e ^\frac{n-2}{2}}-\frac{n-2}{2}\frac{\ln(1+\l_\e^2|x-a_\e|^2)}{\ln \l_\e ^\frac{n-2}{2}}
  +O\left(\frac{e}{\a_\e\d_\e (x)}\frac{1}{\ln \l_\e ^\frac{n-2}{2}}\right)
  \\=&\frac{\ln(\a_\e c_0)}{\ln \l_\e ^\frac{n-2}{2}}-\frac{n-2}{2}\frac{\ln(1+\l_\e^2|x-a_\e|^2)}{\ln \l_\e ^\frac{n-2}{2}}
  +O\left(\frac{1}{\l_\e^\frac{n-2}{4}\ln \l_\e
  ^\frac{n-2}{2}}\right).
 \end{align*}
 The second claim follows from the first one.
  \end{pf}

Now, we need to introduce the following expansion
\begin{lem}\label{lemmaA1} We have
\begin{align}
 & \e \int_{\O}\d_\e^p\ln \left[\frac{\ln(e+\a_\e \delta_\e)}{\ln \l_\e ^\frac{n-2}{2}}\right]\lambda_{\varepsilon} \frac{\partial P \delta_{\varepsilon}}{\partial \lambda}= \Gamma_1 \frac{\e}{\ln \lambda_\e}+o\left(\frac{\e}{\ln\l_\e}\right)+o\left(\frac{1}{(\l_\e d_\e)^{n-2}}\right) , \label{eqa1}\\
 &    \e \int_{\O} \d_\e^p \Big| \ln \left[\frac{\ln(e+\a_\e \delta_\e)}{\ln \l_\e ^\frac{n-2}{2}}\right] \Big|
 \varphi _\e = o\left(\frac{\e}{\ln\l_\e}\right)+o\left(\frac{1}{(\l_\e d_\e)^{n-2}}\right) , \label{eqa2}
\end{align}
where $$\Gamma_1=-\int_{\mathbb{R}^n}\d_\e^p\lambda_{\varepsilon}
\frac{\partial  \delta_{\varepsilon}}{\partial
\lambda}\ln(1+\l_\e^2|x-a_\e|^2)=\displaystyle
c_0\int_{\mathbb{R}^n} \delta_{(0,1)}^p(y)\ln
(\delta_{(0,1)}(y))\frac{1-|y|^2}{(1+|y|^2)^{n/2}}\, dy>0 . $$
\end{lem}
\begin{pf}  We start by proving Assertion \eqref{eqa1}.
 We split the integral over $B_\eta$ and $\O \setminus B_\eta$ where $ B_\eta$ is introduced in Lemma \ref{a34}. On one
hand, using \eqref{lnln} and the fact that $|\lambda_{\varepsilon}
\frac{\partial P \delta_{\varepsilon}}{\partial \lambda}|\leq c
\d_\e$ we have \be \label{eqa4} \displaystyle \int_{\O\setminus
B_\eta}\d_\e^p  \left| \ln \left[\frac{\ln(e+\a_\e \delta_\e)}{\ln
\l_\e ^\frac{n-2}{2}}\right] \right| \left|\lambda_{\varepsilon}
\frac{\partial P \delta_{\varepsilon}}{\partial \lambda}\right| \leq
C \ln\ln \l_\e
 \int_{\O\setminus B_\eta}\d_\e^{p+1}  \leq  C \frac{\ln \ln \l_\e}{(\l_\e\eta_\e)^n}.
\ee
 Now, using the definition of $ \eta_\e $ and using Claim (1) of
Lemma \ref{Lemma 3.3}, we derive that
$$ \frac{\e \ln \ln \l_\e }{ ( \l_\e \eta_\e )^n } \leq \frac{\e \ln \ln \l_\e }{ ( \l_\e d_\e )^n } + \frac{\e \ln \ln \l_\e }
{ ( \l_\e \l_\e^{-3/4} )^n } \leq \frac{ c }{ ( \l_\e d_\e )^n } + \frac{\e \ln \ln \l_\e }{  \l_\e ^{n/4} }
 = o \Big( \frac{ 1 }{ ( \l_\e d_\e )^{n-2} }  + \frac{\e  }{  \ln \l_\e  }  \Big). $$

 Thus \begin{equation}\label{1lem41} \e  \int_{\O\setminus B_\eta}\d_\e^p \left| \ln \left[\frac{\ln(e+\a_\e \delta_\e)}{\ln \l_\e ^\frac{n-2}{2}}\right] \right| \left|\lambda_{\varepsilon} \frac{\partial P \delta_{\varepsilon}}{\partial \lambda}\right|=
 O\left(\frac{\e\ln \ln \l_\e}{(\l_\e\eta_\e)^n} \right)  = o \Big( \frac{ 1 }{ ( \l_\e d_\e )^{n-2} }  + \frac{\e  }{  \ln \l_\e  }  \Big) .
 \end{equation}
 On the other hand, since $B_\eta$ is a subset in $ B_\l=B(a_\e,1/\l_\e^{3/4})$, Eq. \eqref{1ln} holds in
  $B_\eta$.

Furthermore, arguing as in the proof of \eqref{1ln} and using the
fact that $ | \ln t - (t-1) | \leq c | t-1|^2$ in each compact set $
[c_1,c_2] $ with $ c_1>0$, we derive
 \begin{align}\label{ln2}\left|\ln \left[\frac{\ln(e+\a_\e \delta_\e(x))}{\ln \l_\e ^\frac{n-2}{2}}\right]-\left[\frac{\ln(e+\a_\e \delta_\e(x))}{\ln \l_\e ^\frac{n-2}{2}}\right]+1\right|  \leq & c \left|\frac{\ln(e+\a_\e \delta_\e)}{\ln \l_\e ^\frac{n-2}{2}}
  -1\right|^2 \qquad  \forall  \,  x\in B_\eta.
  \end{align}
First, we write
 \begin{align}\label{2lem41}
\displaystyle \int_{B_\eta}\d_\e^p \ln \left[\frac{\ln(e+\a_\e
\delta_\e)}{\ln \l_\e ^\frac{n-2}{2}}\right] \lambda_{\varepsilon}
\frac{\partial P \delta_{\varepsilon}}{\partial \lambda_\e}
=&\int_{B_\eta}\d_\e^p \ln \left[\frac{\ln(e+\a_\e \delta_\e)}{\ln
\l_\e ^\frac{n-2}{2}}\right] \lambda_{\varepsilon} \frac{\partial
\delta_{\varepsilon}}{\partial \lambda_\e} - \int_{B_\eta}\d_\e^p
\ln \left[\frac{\ln(e+\a_\e \delta_\e)}{\ln \l_\e
^\frac{n-2}{2}}\right]\lambda_{\varepsilon}
\frac{\partial \varphi_{\varepsilon}}{\partial \lambda_\e}\nonumber\\
=&I_1 - I_2 .
\end{align}
  Using  \eqref{a4} and Proposition \ref{varphi}, we have
  \be \label{3lem41}
  |I_2|\leq C \Big |\l_\e\frac{\partial \varphi_\e}{\partial
  \l_\e}\Big |_{L^\infty(B_\eta)}
 \frac{1}{\ln \l_\e}\int_{B_\eta}\d_\e^{p}\left[1+\ln \left(1+ \l^2|x-a_\e|^2\right)\right]\,dx
 \leq \frac{C}{(\l_\e d_\e)^{n-2}\ln\l_\e}.
  \ee
  To compute $I_1$, expanding $ \ln (1+t) = t + O(t^2) $ in  $[-c_1, c_2]$ (with $c_1 < 1$) and using Lemma \ref{a34}, we get
  \begin{align}
  I_1=&\int_{B_\eta}\d_\e^p \ln \left[1+\left(\frac{\ln(e+\a_\e \delta_\e)}{\ln \l_\e ^\frac{n-2}{2}}-1\right)\right]\lambda_{\varepsilon} \frac{\partial  \delta_{\varepsilon}}{\partial \lambda}\nonumber\\
  =&\int_{B_\eta}\d_\e^p \left(\frac{\ln(e+\a_\e \delta_\e)}{\ln \l_\e ^\frac{n-2}{2}}-1\right)\lambda_{\varepsilon} \frac{\partial  \delta_{\varepsilon}}{\partial \lambda}+O\left(\int_{B_\eta}\d_\e^p \left(\frac{\ln(e+\a_\e \delta_\e)}{\ln \l_\e ^\frac{n-2}{2}}-1\right)^2\lambda_{\varepsilon} \frac{\partial  \delta_{\varepsilon}}{\partial \lambda}\right)\nonumber\\
  =&\frac{\ln(\a_\e c_0)}{\ln \l_\e ^\frac{n-2}{2}}\int_{B_\eta}\d_\e^p\lambda_{\varepsilon} \frac{\partial  \delta_{\varepsilon}}{\partial \lambda}
  -\frac{n-2}{2}\frac{1}{\ln \l_\e ^\frac{n-2}{2}}\int_{B_\eta}
  \d_\e^p\lambda_{\varepsilon} \frac{\partial  \delta_{\varepsilon}}{\partial \lambda}\ln(1+\l_\e^2|x-a_\e|^2)
  \nonumber\\
  & + O\left(\frac{1}{\l_\e^\frac{n-2}{4}\ln \l_\e
  ^\frac{n-2}{2}}\int_{B_\eta}\d_\e^{ p +1 } \right)+O\left(\int_{B_\eta}\d_\e^{p+1}
  \left[\frac{1+\ln(1+\l_\e^2|x-a_\e|^2)}{\ln \l_\e ^\frac{n-2}{2}}\right]^2\right)\nonumber\\
 =&I_{11}+I_{12}+I_{13}+I_{14}.
  \end{align}
   In the following, we compute each integral.
  \begin{align}
  I_{11} & =\frac{\ln(\a_\e c_0)}{\ln \l_\e ^\frac{n-2}{2}}\left(\int_{\mathbb{R}^n}\ldots -
  \int_{B_\eta^c}\ldots \right)=O\left(\frac{1}{\ln\l_\e (\l_\e\eta_\e)^n}\right) , \\
  I_{12}  & =  - \frac{1}{\ln\l_\e}\left(\int_{\mathbb{R}^n}\d_\e^p\lambda_{\varepsilon} \frac{\partial  \delta_{\varepsilon}}{\partial \lambda}\ln(1+\l_\e^2|x-a_\e|^2) - \int_{B_\eta^c}\d_\e^p\lambda_{\varepsilon} \frac{\partial  \delta_{\varepsilon}}{\partial \lambda}\ln(1+\l_\e^2|x-a_\e|^2)\right)\nonumber\\
   & =  \frac{\Gamma_1}{\ln\l_\e}+O\left(\frac{1}{\ln\l_\e (\l_\e\eta_\e)^{n-1}}\right)
  \end{align}
  where $\Gamma_1$ is the positive constant introduced in Lemma  \ref{lemmaA1}. In addition, we have
  \begin{equation}
  I_{13}=O\left(\frac{1}{\l_\e^\frac{n-2}{4}\ln\l_\e}\right) \qquad \mbox{ and } \qquad
  I_{14}=O\left(\frac{1}{(\ln\l_\e)^2}\right).\label{6lem41}
  \end{equation}
Combining \eqref{1lem41} and \eqref{2lem41}-\eqref{6lem41}, the
proof of Assertion \eqref{eqa1} of Lemma \ref{lemmaA1} follows.

 Concerning the proof of Assertion \eqref{eqa2}, it can be deduced from the proof of equations \eqref{eqa4}, \eqref{1lem41} and \eqref{3lem41}. \\
Hence the proof of the lemma is completed.
\end{pf}
\begin{pro}\label{Proposition 4.7} Let $u_{\varepsilon}$ satisfy the assumption of Theorem 1.4. Then there exist $C_{1}>0$ and $C_{2}>0$ such that
$$
C_1 \frac{\varepsilon}{\ln \l_\e}(1+o(1))+ C_2
\frac{H\left(a_{\varepsilon},
a_{\varepsilon}\right)}{\lambda_{\varepsilon}^{n-2}}(1+o(1)) = o
\left(\frac{1}{(\l_\e d_\e)^{n-2} }\right) .
$$
\end{pro}
\begin{pf}  Multiplying the equation $\left(SC_\e\right)$
by $\lambda_{\varepsilon}\left(\partial P
\delta_{\varepsilon}\right) /\left(\partial
\lambda_{\varepsilon}\right)$ and integrating over $\Omega$, we
obtain
\begin{align}\label{2.10}
0= & \int_{\Omega} (-\Delta) u_{\varepsilon} \lambda_{\varepsilon} \frac{\partial P \delta_{\varepsilon}}{\partial \lambda_\e} - \int_{\Omega} | u_{\varepsilon} | ^{p-1} u_\e \ln (e+ | u_\e | )^\e \lambda_{\varepsilon} \frac{\partial P \delta_{\varepsilon}}{\partial \lambda_\e}\nonumber \\
= & \alpha_{\varepsilon} \int_{\Omega} \delta_{\varepsilon}^{p}
\lambda_{\varepsilon} \frac{\partial P
\delta_{\varepsilon}}{\partial \lambda_\e}  - \int_{\Omega} | u_\e
|^{p-1} u_\e g_\e(u_\e)\l_\e\frac{\partial P
\delta_{\varepsilon}}{\partial \lambda_\e}
\end{align}
since $v_\e \in E_{\left(a_{\varepsilon},\l_\e\right)}$.\\
We estimate each term on the right-hand side in \eqref{2.10}. First,
from \cite[Lemma 2.2]{BFG}, we have
\begin{equation}\label{dev1}
 \langle  P \delta_{\varepsilon} ,
\lambda_{\varepsilon} \frac{\partial P
\delta_{\varepsilon}}{\partial \lambda_\e} \rangle =
 \int_{\Omega} \delta_{\varepsilon}^{p} \lambda_{\varepsilon} \frac{\partial P \delta_{\varepsilon}}{\partial \lambda_\e}=
 \frac{n-2}{2}c_1\frac{H(a_\e,a_\e)}{\l_\e^{n-2}}+O\left(\frac{\ln(\l_\e d_\e)}{(\l_\e d_\e)^n}\right)
 \end{equation}
with $c_{1}=c_{0}^{2 n /(n-2)} \int_{\mathbb{R}^{n}} \frac{d x}{\left(1+|x|^{2}\right)^{(n+2) / 2}}$.\\
 Second, we write
\begin{align}\label{2.12}
 \int_{\Omega} | u_\e |^{p-1} u_\e g_\e(u_\e) \lambda_{\varepsilon} \frac{\partial P \delta_{\varepsilon}}{\partial \lambda_\e} & =
 \int_{\Omega} | u_\e |^{p-1} u_\e \left[g_\e(u_\e) - g_\e(\a_\e P\delta_\e) -g'_\e(\a_\e P\d_\e)v_\e \right]
 \lambda_{\varepsilon} \frac{\partial P \delta_{\varepsilon}}{\partial \lambda_\e} \nonumber \\
 & \quad  +\int_{\Omega} | u_\e |^{p-1} u_\e  g_\e(\a_\e P\delta_\e) \lambda_{\varepsilon}
\frac{\partial P \delta_{\varepsilon}}{\partial \lambda_\e}
+\int_{\Omega} | u_\e |^{p-1} u_\e g'_\e(\a_\e P\d_\e)v_\e \lambda_{\varepsilon} \frac{\partial P \delta_{\varepsilon}}{\partial \lambda_\e} \nonumber\\
& : = A+B+C.
\end{align}
and we have to estimate each term on the right hand-side of \eqref{2.12}.\\
We claim that
\begin{equation}\label{estA}
A=o\left(\|v_\e\|^2\right).
\end{equation}
In fact, let $ \O_v $ be defined in \eqref{Ov1}. We split the
integral $A$ into integrals over $\Omega_v$ and $\Omega_v^c$. Using
\eqref{gB.1}, Sobolev embedding theorem and the fact \be
\label{ldpdless} |\lambda_{\varepsilon} {\partial P
\delta_{\varepsilon}} / {\partial \lambda _\e }|\leq c P\d_\e, \ee
  we have
\begin{align}\label{1A0}
\int_{\Omega_v} | u_\e | ^{p-1} u_\e  g'_\e(\a_\e P\d_\e) v_\e
 \lambda_{\varepsilon} \frac{\partial P \delta_{\varepsilon}}{\partial
 \lambda_\e} & = O\left(\e\int_{\Omega_v} | v_\e | ^{p+1}\frac{ P\d_\e }{e+\a_\e P\d_\e}
   \right)
  = O\left(\e\int_{\Omega_v} | v_\e | ^{p+1}
\right)
  = O(\e\|v\|^{p+1}).
\end{align}
In addition, using  Lemma \ref{integralv}, and arguing as in
\eqref{vest1} (by using \eqref{X3} and, in the same way, the fact
that $ \ln (e + \a_\e P \d_\e) \leq c \ln (e+ | v _\e |)$ in $
\O_v$), we get
\begin{align}\label{A1}
 \int_{\Omega_v} | u_\e | ^{p-1} u_\e \left[g_\e(u_\e) -g_\e(\a_\e P\delta_\e)  \right]
 \lambda_{\varepsilon} \frac{\partial P \delta_{\varepsilon}}{\partial \lambda_\e} & = O\left(\int_{\Omega_v} | u_\e | ^p g_\e(u_\e)| v_\e|+ \int_{\Omega_v} | u_\e | ^p
g_\e(\a_\e P\delta_\e) |v_\e|\right)
  = o( \| v_\e \| ^2 ).
\end{align}

Concerning the integral over $\O\setminus \O_v$, we use the mean
value theorem, \eqref{gB''}, \eqref{vlesspdelta}, \eqref{ldpdless},
Holder inequality and Sobolev embedding theorem, we obtain
\begin{align}\label{A3}
  \int_{\Omega\setminus \O_v} | u_\e |^{p-1} u_\e \left[g_\e(u_\e) - g_\e(\a_\e P\delta_\e) -g'_\e(\a_\e P\d_\e)v_\e \right]
 \lambda_{\varepsilon} \frac{\partial P \delta_{\varepsilon}}{\partial \lambda_\e} & = O\left(\e\int_{\O\setminus \O_v} | u_\e | ^p\frac{ v_\e^2 P\d_\e}{(e+\a_\e P\d_\e+\theta
 v_\e)^2} \right)\nonumber \\
 & = O\left(\e\int_{\O\setminus \O_v} | u_\e | ^{p-1} {v_\e^2}  \right)\nonumber \\
 & = O(\e\|v\|^2).
\end{align}

 Combining \eqref{1A0}-\eqref{A3}, Claim \eqref{estA} follows.

Now, we compute the integral $B$. We split it as follows
\begin{align}\label{B}
 B & = \int_{\Omega} | u_\e | ^{p-1} u_\e  \left[g_\e(\a_\e P\d_\e) -g_\e(\a_\e \d_\e)\right]\lambda_{\varepsilon} \frac{\partial P \delta_{\varepsilon}}{\partial \lambda_\e} + \int_{\Omega} | u_\e | ^{p-1} u_\e  g_\e(\a_\e \delta_\e) \lambda_{\varepsilon} \frac{\partial P \delta_{\varepsilon}}{\partial \lambda_\e} .
\end{align}
By using the mean value theorem, \eqref{3.1}, \eqref{gB.1},
\eqref{pdless}, \eqref{ldpdless}, Proposition \ref{varphi} and Lemma
\ref{Lemma 4.6},  we get
 \begin{align}\label{b1}
 \int_{\Omega} | u_\e | ^{p-1} u_\e  \left[g_\e(\a_\e P\d_\e) -g_\e(\a_\e \d_\e)\right]\lambda_{\varepsilon} \frac{\partial P \delta_{\varepsilon}}{\partial \lambda_\e}  & = O\left(\e\int_{\O} | u_\e | ^p\frac{P\d_\e}{e+\a_\e\d_\e-\theta\a_\e\varphi_\e}\varphi_\e\right)\nonumber  \\
 & = O\left(\e\int_{\O} | u_\e | ^p\varphi_\e\right)\nonumber\\
  & = O\left(\e | \varphi_\e  |_{ L^\infty} \int_{ \O } \d_\e^p + \e\int_{\O}|v_\e|^p\varphi_\e\right)\nonumber\\
 & = O\left(\e\frac{1}{(\l_\e d_\e)^{n-2}}+\frac{\e\|v\|^p}{(\l_\e d_\e)^\frac{n-2}{2}}\right)\nonumber\\
 & = o\Big( \frac{\e}{\ln \l_\e} + \frac{1}{(\l_\e d_\e)^{n-2}} \Big).
 \end{align}

 Now, we will focus on the second integral of \eqref{B}. Let $ \O_1$ and $ \O_2$ be defined in \eqref{O1O2}, we have
    \begin{align}\label{b2}
  \int_{\Omega} | u_\e | ^{p-1} u_\e  g_\e(\a_\e \delta_\e) \lambda_{\varepsilon} \frac{\partial P \delta_{\varepsilon}}{\partial \lambda_\e}  & = \a_\e^p\int_{\Omega}\d_\e^p g_\e\left(\a_\e\delta_{\varepsilon}\right)\lambda_{\varepsilon} \frac{\partial P \delta_{\varepsilon}}{\partial \lambda_\e} + p \a_\e^{p-1}\int_{\Omega}\d_\e^{p-1} ( - \a_\e \varphi_\e + v_\e) g_\e (\a_\e\delta_{\varepsilon} )\lambda_{\varepsilon} \frac{\partial P \delta_{\varepsilon}}{\partial \lambda_\e} \nonumber \\
 & \quad + O\left(\int_{ \O_1} \d_\e^{p-1}  | - \a_\e \varphi_\e + v_\e |  ^2 g_\e\left(\a_\e\delta_{\varepsilon}\right) + \int_{ \O_2 } | -\a_\e \varphi_\e + v_\e | ^p \d_\e g_\e\left(\a_\e\delta_{\varepsilon}\right)\right) \nonumber \\
 & = B_{1} + B_{2} + O (  B_3 + B_4 ).
\end{align}
Observe that, using \eqref{majln} and \eqref{3.2}, we deduce that $
g_\e ( \a_\e \d_\e) $ is bounded. Thus, using Proposition
\ref{varphi} and Lemma \ref{Lemma 4.6}, we get
\begin{align}
 & B_4 \leq c \int_\O ( \varphi _\e ^{p+1} + | v _\e |^{ p+1}) \leq c \| \varphi _\e \|_{ L^{p+1} } ^{p+1} + c  \| v _\e \|^{p+1} = o \Big( \frac{ \e }{ \ln \l_\e} +  \frac{1}{( \l_\e d_\e)^{n-2} } \Big) , \label{NB4}  \\
 & B_3 \leq c \int_{ B(a_\e , d_\e)} \d_\e^{p-1}   \varphi_\e ^2  + \int_{ \O \setminus B(a_\e , d_\e)} \d_\e ^{ p+1} + \int_\O \d_\e^{p-1}  v_\e  ^2 \notag  \\
  & \quad \leq | \varphi_\e |_{L^\infty} ^2  \int_{ B(a_\e , d_\e)} \d_\e^{p-1} + \frac{c } { ( \l_\e d_\e )^n } + c \| v_\e \|^2  \notag \\
   & \quad  = o \Big( \frac{ \e }{ \ln \l_\e} +  \frac{1}{( \l_\e d_\e)^{n-2} } \Big) . \label{NB3}
\end{align}
Concerning $ B_2$, we split it into two pieces. The first one
contains the $ v_\e$.  Using Lemmas \ref{Lemma 3.3}, \ref{Lemma
4.6}, \eqref{2V}, Proposition \ref{varphi} and the fact that $ v_\e
\in E_{(a_\e , \l_\e)}$, we have
\begin{align}\label{b22}
 \int_{\Omega} \d_\e^{p-1}  v_\e g_\e (\a_\e\delta_{\varepsilon} )\lambda_{\varepsilon} \frac{\partial P \delta_{\varepsilon}}{\partial \lambda_\e}& = \ln ^\e (\l_\e^\frac{n-2}{2})\int_{\Omega}\delta_{\varepsilon}^{p-1}\lambda_{\varepsilon} \frac{\partial  P\delta_{\varepsilon}}{\partial \lambda _\e }v_\e+O\left(\e\int_{\O}\d_\e^p|v_\e| \Big| \ln \left[\frac{\ln(e+\a_\e \delta_\e)}{\ln \l_\e ^\frac{n-2}{2}}\right] \Big| \right)\nonumber \\
 & = \ln ^\e (\l_\e^\frac{n-2}{2})\int_{\Omega}\delta_{\varepsilon}^{p-1}\lambda_{\varepsilon} \frac{\partial  \delta_{\varepsilon}}{\partial \lambda _\e }v_\e + O\Big( \int_\O \d_\e^{p-1} \l_\e \Big| \frac{ \partial \varphi _\e}{ \partial \l_\e}\Big|  | v_\e | + \|v_\e\|\frac{\e}{\ln\l_\e}\Big) \notag \\
  & = O\Big( \| v_\e\| \Big| \l_\e \frac{ \partial \varphi _\e}{ \partial \l_\e}\Big| _{L^{p+1}} + \|v_\e\|\frac{\e}{\ln\l_\e}\Big) \notag \\
   &  = o \Big( \frac{ \e }{ \ln \l_\e} +  \frac{1}{( \l_\e d_\e)^{n-2} } \Big)  .
\end{align}
Regarding the second part of $ B_2$, using Lemmas \ref{Lemma 3.3},
we obtain
\begin{align} \label{wx3}
 & p  \int_{\Omega} \d_\e^{p-1}  \varphi_\e g_\e (\a_\e\delta_{\varepsilon} )\lambda_{\varepsilon} \frac{\partial P \delta_{\varepsilon}}{\partial \lambda_\e} \notag \\
  & = \ln ^\e (\l_\e^\frac{n-2}{2})  p \int_{\Omega} \d_\e^{p-1}  \varphi_\e \lambda_{\varepsilon} \frac{\partial  \delta_{\varepsilon}}{\partial \lambda_\e} + O \Big( \e \int_{\Omega} \d_\e^{p}  \varphi_\e \Big| \ln \left[\frac{\ln(e+\a_\e \delta_\e)}{\ln \l_\e ^\frac{n-2}{2}}\right] \Big| +  \int_{\Omega} \d_\e^{p-1}  \varphi_\e\Big|  \lambda_{\varepsilon} \frac{\partial  \varphi_{\varepsilon}}{\partial \lambda_\e} \Big| \Big)
 \end{align}
Observe that
$$ \int_\O \d_\e ^{ p} \l_\e \frac{ \partial \d_\e}{ \partial \l_\e} = O \Big( \frac{1} { ( \l_\e d_\e)^n} \Big) $$
which implies that \be \label{wx1} p \int_{\Omega} \d_\e^{p-1}
\varphi_\e \lambda_{\varepsilon} \frac{\partial
\delta_{\varepsilon}}{\partial \lambda_\e}  = - \langle
\lambda_{\varepsilon} \frac{\partial P
\delta_{\varepsilon}}{\partial \lambda_\e} , P\d_\e \rangle + O
\Big( \frac{1} { ( \l_\e d_\e)^n} \Big) . \ee In addition, using
Proposition \ref{varphi}, we have \be \label{wx2} \int_{\Omega}
\d_\e^{p-1}  \varphi_\e\Big|  \lambda_{\varepsilon} \frac{\partial
\varphi_{\varepsilon}}{\partial \lambda_\e} \Big| \leq |\varphi_\e
|_{L^\infty} \Big| \lambda_{\varepsilon} \frac{\partial
\varphi_{\varepsilon}}{\partial \lambda_\e} \Big|_{L^\infty}
\int_{B(a_\e , d_\e) }  \d_\e ^{p-1}  + \int_{ \O \setminus B(a_\e ,
d_\e)} \d_\e ^{p+1} = o\Big( \frac{\e }{ \ln \l_\e} + \frac{ 1 }{ (
\l_e d_\e )^{n-2}} \Big) . \ee Thus, combining \eqref{eqa2},
\eqref{dev1} and \eqref{b22}-\eqref{wx2}, the estimate of $ B_2$
becomes \be \label{NB2} B_2 = \a_\e^p  \ln ^\e (\l_\e^\frac{n-2}{2})
\frac{n-2}{2}c_1\frac{H(a_\e,a_\e)}{\l_\e^{n-2}}  +  o\Big( \frac{\e
}{ \ln \l_\e} + \frac{ 1 }{ ( \l_e d_\e )^{n-2}} \Big)  . \ee

Now, we will focus on estimating the integral $ B_1$.  Using Lemma
\ref{Lemma 3.3}, we write
\begin{align}\label{b21}
B_{1} & =
\a_\e^p\int_{\Omega}\d_\e^pg_\e\left(\a_\e\delta_{\varepsilon}\right)\lambda_{\varepsilon}
\frac{\partial P \delta_{\varepsilon}}{\partial \lambda_\e}\nonumber\\
 & =  \a_\e^p\ln ^\e
(\l_\e^\frac{n-2}{2})
\int_{\Omega}\delta_{\varepsilon}^p\lambda_{\varepsilon}
\frac{\partial  P\delta_{\varepsilon}}{\partial \lambda_\e}+
\a_\e^p\e\int_{\O}\d_\e^p\ln \left[\frac{\ln(e+\a_\e \delta_\e)}{\ln
\l_\e ^\frac{n-2}{2}}\right] \lambda_{\varepsilon} \frac{\partial P
\delta_{\varepsilon}}{\partial
\lambda_\e}+O\left(\e^2\int_{\O}\d_\e^{p+1}
\ln \left[\frac{\ln(e+\a_\e \delta_\e)}{\ln \l_\e ^\frac{n-2}{2}}\right]^2\right) \nonumber\\
 & := B_{11}+B_{12}+ O( B_{13}).
\end{align}
Through \eqref{dev1} and \eqref{Lemma 4.2}, we get
\begin{align}\label{b211}
B_{11}=&\a_\e\frac{n-2}{2}c_1\frac{H(a_\e,a_\e)}{\l_\e^{n-2}}(1+o(1))+O\left(\frac{\ln(\l_\e
d_\e)}{(\l_\e d_\e)^n}\right).
\end{align}
    Using Lemma \ref{lemmaA1}, we have
    \begin{equation}\label{b212}B_{12}=  \a_\e^p\, \Gamma_1 \frac{\e}{\ln \lambda_\e}  + o\left(\frac{\e}{\ln\l_\e} + \frac{1}{(\l_\e d_\e)^{n-2}}\right). \end{equation}
To estimate $B_{13}$, we recall the set $B_\eta =B(a_\e,\eta_\e)=
B(a_\e,\l_\e^{-3/4})\cap B(a_\e,d_\e)$. Using \eqref{lnln}, we have
\be\label{b2131} \e^2 \int_{\O\setminus B_\eta}\d_\e^{p+1} \ln
\left[\frac{\ln(e+\a_\e \delta_\e)}{\ln \l_\e
^\frac{n-2}{2}}\right]^2 \leq C \, \e^2  \ln(\ln \l_\e)^2
 \int_{\O\setminus B_\eta}\d_\e^{p+1}  \leq  C \, \e^2 \frac{\ln (\ln \l_\e)^2}{(\l_\e\eta_\e)^n}  = o\left(\frac{\e}{\ln\l_\e} + \frac{1}{(\l_\e d_\e)^{n-2}}\right) .
\ee Note that $ B_\eta \subset B_\l := B(a_\e, \l_\e ^{-3/4})$, and
using \eqref{1ln}, we have
  \begin{align}\label{b2132}
  \e^2\displaystyle \int_{ B_\l}\d_\e^{p+1} \ln \left[\frac{\ln(e+\a_\e \delta_\e)}{\ln \l_\e ^\frac{n-2}{2}}\right]^2\leq&
C \left(\frac{\e}{\ln
\l_\e}\right)^2\left(\int_{B_\l}\d_\e^{p+1}\left|\ln
\left(\frac{e}{\l_\e^\frac{n-2}{2}}+
 \frac{\a_\e c_0}{(1+\l_\e^2|x-a_\e|^2)^\frac{n-2}{2}}\right)\right|^2\, dx\right)\nonumber\\
\leq& C \left(\frac{\e}{\ln
\l_\e}\right)^2\int_{\tilde{B}_\l}\frac{1}{(1+|y|^2)^n}\left|\ln
\left(\frac{e}{\l_\e^\frac{n-2}{2}}+ \frac{\a_\e
c_0}{(1+|y|^2)^\frac{n-2}{2}}\right)\right|^2\,
dy\nonumber\\
\leq& C \left(\frac{\e}{\ln \l_\e}\right)^2
  \end{align}
 by arguing as in the proof of \eqref{5V}, where $\tilde{B}_\l=B(0, \l_\e ^{1/4})$ and we have used the change of coordinates  $y=\l_\e(x-a_\e)$. Hence \eqref{b2131} and \eqref{b2132} assert  that
  \begin{equation}\label{b213}
  B_{13} =  o\left(\frac{\e}{\ln\l_\e} + \frac{1}{(\l_\e d_\e)^{n-2}}\right).
   \end{equation}
 Equations \eqref{B}- \eqref{NB3}, \eqref{NB2}-
\eqref{b212} and \eqref{b213} imply that
\begin{equation}\label{estB} B =  \a_\e^p\, \Gamma_1\displaystyle
\frac{\e}{\ln \lambda_\e}  + o\left(\frac{\e}{\ln\l_\e}\right)+ (
\a_\e +  \a_\e^p  \ln ^\e (\l_\e^\frac{n-2}{2})  )
\frac{n-2}{2}c_1\frac{H(a_\e,a_\e)}{\l_\e^{n-2}}(1+o(1))+o\left(\frac{1}{(\l_\e
d_\e)^{n-2}}\right).
\end{equation}

Lastly, we estimate $C$. We start by splitting it as follows
\begin{align}\label{2.13}
C & = \int_{\Omega} | u_\e | ^{p-1} u_\e g'_\e(\a_\e P\delta_\e)
v_\e
\lambda_{\varepsilon}\frac{\partial P \delta_{\varepsilon}}{\partial \lambda} \nonumber \\
 & = \int_{\Omega} | u_\e | ^{p-1} u_\e \left[g'_\e(\a_\e P\d_\e) -g'_\e(\a_\e \d_\e)\right]v_\e \lambda_{\varepsilon} \frac{\partial P \delta_{\varepsilon}}
{\partial \lambda} + \int_{\Omega} | u_\e | ^{p-1} u_\e g'_\e(\a_\e \delta_\e) v_\e \lambda_{\varepsilon} \frac{\partial P \delta_{\varepsilon}}{\partial \lambda}\nonumber\\
 & := C_1+C_2.
\end{align}
We claim that
\begin{equation}\label{estC}
C=O\left(\e\|v_\e\|^2+\|v_\e\|\left(\frac{\e}{\ln
\l_\e}+\frac{1}{(\l_\e d_\e)^{n-2}}\right)\right). \end{equation}
Indeed, by the mean value theorem, \eqref{gB.1}, \eqref{gB''},
\eqref{ldpdless}, \eqref{2vest2'}, Holder inequality, Sobolev
embedding theorem and Lemma \ref{Lemma 4.6}, we have
\begin{align}\label{estC1}
 C_1 & = O\left(\e\int_{\Omega}P\d_\e^{p+1}\frac{1}{(e+\a_\e\d_\e-\theta \varphi_\e)^2}\varphi_\e |v_\e |  +\e\int_{\Omega}|v_\e|^{p+1}\frac{1}{e+\a_\e P\d_\e}  P \delta_{\varepsilon}\right)\nonumber \\
 & = O\left(\e\int_{\Omega}\d_\e^{p-1}\varphi_\e |v_\e |  +\e\int_{\Omega}|v_\e|^{p+1}\right)\nonumber \\
 & = O\left(\e\int_{B_\e}\d_\e^{p-1}\varphi_\e| v_\e|+\e\int_{B_\e ^c }\d_\e^{p} v_\e  +\e\int_{\Omega}|v_\e|^{p+1}\right)\nonumber  \\
  &  = o\Big( \frac{\e }{ \ln \l_\e } + \frac{ 1 }{ ( \l_\e d_\e )^{n-2} } \Big)
\end{align}
where $ B_\e := B( a_\e , d_\e/2) $. \\
Recall that $\Omega_\l:= B(a_\e, \l_\e^{-3/4})$. Taking account of
\eqref{ldpdless} and \eqref{gB.1}, we write
\begin{align}\label{estC2}
C_2 & = O\left(\int_{\Omega}P\d_\e^{p+1} |g'_\e (\a_\e\delta_\e)| |v_\e|+\e\int_{\Omega}|v_\e|^{p+1} \right)\nonumber\\
 & = O\left(\int_{\Omega_\l}\d_\e^{p+1} |g'_\e (\a_\e\delta_\e)| |v_\e|+\int_{\Omega\setminus\Omega_\l}\d_\e^{p+1} |g'_\e
(\a_\e\delta_\e)| |v_\e|+\e\|v_\e\|^{p+1} \right).
\end{align}
 On one hand, we recall that in $\O_\l$, we have $\delta_{\varepsilon}\geq c\l_\e^{(n-2)/4}$. Thus
 \begin{equation}\label{0estC}
\ln(e+\a_\e \d_\e)\geq c\ln \l_\e \hbox{ in }\O_\l.
 \end{equation}Using \eqref{g'}, \eqref{0estC} and \eqref{3.2}, Holder
inequality and Sobolev embedding theorem, we get
\begin{align}\label{estC21}
\int_{\Omega_\l}\d_\e^{p+1} |g'_\e (\a_\e\delta_\e)| |v_\e| & =
O\left(\e\int_{\Omega_\l}\d_\e^{p+1} \frac{\ln(e+\a_\e
\d_\e)^{\e-1}}{e+\a_\e \d_\e} |v_\e|\right)
\nonumber\\
 & = O\left(\e\ln(\l_\e)^{\e-1}\int_{\Omega_\l}\d_\e^{p} |v_\e|\right) \nonumber\\
 & = O\left(\frac{\e}{\ln \l_\e}\|v_\e\|\right) .
\end{align}
On the other hand, taking account of \eqref{gB.1}, Holder inequality
and Sobolev embedding theorem, we obtain
\begin{align}\label{estC22}
\int_{\O\setminus\Omega_\l}\d_\e^{p+1} |g'_\e (\a_\e\delta_\e)|
|v_\e|=&O\left(\e\int_{\O\setminus\Omega_\l}\d_\e^{p+1}
\frac{1}{e+\a_\e \d_\e} |v_\e|\right)
\nonumber\\
=&O\left(\e\int_{\O\setminus\Omega_\l}\d_\e^{p} |v_\e|\right)
\nonumber\\
=&O\left(\frac{\e}{ \l_\e^\frac{n+2}{8}}\|v_\e\|\right).
\end{align}
Thus we derive \eqref{estC} from combining \eqref{2.13}, \eqref{estC1}, \eqref{estC2}, \eqref{estC21}, \eqref{estC22} and Lemma \ref{Lemma 4.6}. \\
 \eqref{2.10} together with \eqref{dev1},
\eqref{2.12}, \eqref{estA}, \eqref{estB}, \eqref{estC} and the fact
that $\a_\e$ satisfies \eqref{3.1} and Lemma \ref{Lemma 4.2}, the
result of Proposition \ref{Proposition 4.7} follows.
\end{pf}


We are now able to prove Theorem \ref{1.4}.

{\bf{Proof of Theorem \ref{1.4}}} Arguing by contradiction, let us
suppose that $\left(SC_\e\right)$ has a solution $u_{\varepsilon}$
as stated in Theorem \ref{1.4}. From Proposition \ref{Proposition
4.7}, we have

\begin{align}\label{(4.20)}
C_1 \frac{\varepsilon}{\ln \l_\e}(1+o(1))+C_2
\frac{H\left(a_{\varepsilon},
a_{\varepsilon}\right)}{\lambda_{\varepsilon}^{n-2}}(1+o(1))=o\left(\frac{1}{(\l_\e
d_\e)^{n-2}}\right) \end{align}
with $C_{1}>0$ and $C_{2}>0$.\\
Two cases may occur :

Case $1$. $d_{\varepsilon} \rightarrow 0$ as $\varepsilon
\rightarrow 0$. Using \eqref{(4.20)} and the fact that
$H\left(a_{\varepsilon}, a_{\varepsilon}\right) = c
d_{\varepsilon}^{2-n} (1+o(1)) $, we derive a contradiction.

Case $2$. $d_{\varepsilon} \nrightarrow 0$ as $\varepsilon
\rightarrow 0$.
We have $H\left(a_{\varepsilon}, a_{\varepsilon}\right) \geq c>0$ as $\varepsilon \rightarrow 0$ and \eqref{(4.20)} also leads to a contradiction.\\
 This completes the proof of Theorem \ref{1.4}.



\end{document}